\documentclass[letterpaper,leqno,11pt,oneside]{amsart}

\usepackage{amsmath,amsthm,amsfonts,amssymb} \usepackage{enumitem}
\usepackage[usenames]{color} \usepackage{mathtools}
\usepackage[extension=pdf]{hyperref} \usepackage[height=9.6in,width=5.95in]{geometry}
\usepackage{graphics,graphicx}
\usepackage[labelfont=rm]{subcaption}
\usepackage{xr}

\usepackage[style=alphabetic,natbib=true,sorting=nyt,maxnames=99,isbn=false,doi=true,url=false,firstinits=true,hyperref=auto,arxiv=abs,backend=bibtex]{biblatex}
\DeclareFieldFormat[article,inbook,incollection,inproceedings,patent,thesis,unpublished]{title}{#1\isdot}
\renewbibmacro{in:}{%
  \ifentrytype{article}{}{%
    \printtext{\bibstring{in}\intitlepunct}}} 
\defbibheading{apa}[\refname]{\section*{#1}}
\DeclareFieldFormat{sentencecase}{\MakeSentenceCase{#1}} \renewbibmacro*{title}{%
  \ifthenelse{\iffieldundef{title}\AND\iffieldundef{subtitle}}{}
  {\ifthenelse{\ifentrytype{article}\OR\ifentrytype{inbook}%
      \OR\ifentrytype{incollection}\OR\ifentrytype{inproceedings}%
      \OR\ifentrytype{inreference}} {\printtext[title]{%
        \printfield[sentencecase]{title}%
        \setunit{\subtitlepunct}%
        \printfield[sentencecase]{subtitle}}}%
    {\printtext[title]{%
        \printfield[titlecase]{title}%
        \setunit{\subtitlepunct}%
        \printfield[titlecase]{subtitle}}}%
    \newunit}%
  \printfield{titleaddon}}

\definecolor{darkblue}{rgb}{0.13,0.13,0.39}
\hypersetup{colorlinks=true,urlcolor=darkblue,citecolor=darkblue,linkcolor=darkblue,%
  pdftitle=Airy processes and variational problems}

\mathtoolsset{showonlyrefs,showmanualtags}

\newtheorem{thm}{Theorem}[section] \newtheorem{lem}[thm]{Lemma}

\theoremstyle{definition} \newtheorem{rem}[thm]{Remark} \newtheorem*{rem*}{Remark}
 
\newtheorem*{exs}{Examples}
 \newcounter{assum}

\newcommand{\I}{{\rm i}} \newcommand{\pp}{\mathbb{P}} 
 \newcommand{\ee}{\mathbb{E}} \newcommand{\rr}{\mathbb{R}}
\newcommand{\nn}{\mathbb{N}} \newcommand{\zz}{\mathbb{Z}} \newcommand{\aip}{\mathcal{A}_2}
\newcommand{\aipo}{\mathcal{A}_1} \newcommand{\Bt}{{\mathcal{A}}_{2\to 1}}
  \newcommand{\ct}{\mathcal{T}}
\newcommand{\cx}{\mathcal{X}}
\newcommand{\cm}{\mathcal{M}}  \newcommand{\p}{\partial}
\newcommand{\uno}[1]{\mathbf{1}_{#1}}

\newcommand{\ep}{\varepsilon}  \newcommand{\wt}{\widetilde}
 
\newcommand{\K}{K_{\Ai}} \newcommand{\m}{\overline{m}}

 \newcommand{\qand}{\quad\text{and}\quad}
\newcommand{\qqand}{\qquad\text{and}\qquad}
\DeclareMathOperator{\Ai}{Ai} \DeclareMathOperator{\tr}{tr}
\newcommand{\e}{\epsilon}
\newcommand{\z}{\zeta}

\newcommand{\gref}[1]{\ref*{g-#1} of \cite{cqr}}  \newcommand{\grefn}[1]{\ref*{g-#1}} 
\newcommand{\pref}[1]{\ref*{p-#1} of \cite{mqr}} \newcommand{\peqref}[1]{(\ref*{p-#1}) in
  \cite{mqr}} \newcommand{\prefn}[1]{\ref*{p-#1}} 
  \newcommand{\orefn}[1]{\ref*{o-#1}}

\DeclareMathOperator*{\argmax}{arg\,max}

\numberwithin{equation}{section}

\let\oldmarginpar\marginpar
\renewcommand\marginpar[1]{\-\oldmarginpar[\raggedleft\footnotesize #1]%
  {\raggedright{\small\textsf{#1}}}}

\allowdisplaybreaks[2]

\begin{document}

\title{Airy Processes and Variational Problems}

\author{Jeremy Quastel} \address[J.~Quastel]{
  Department of Mathematics\\
  University of Toronto\\
  40 St. George Street\\
  Toronto, Ontario\\
  Canada M5S 2E4} \email{quastel@math.toronto.edu} \author{Daniel Remenik}
\address[D.~Remenik]{
  Departamento de Ingenier\'ia Matem\'atica and Centro de Modelamiento Matem\'atico\\
  Universidad de Chile\\
  Av. Blanco Encala\-da 2120\\
  Santiago\\
  Chile} \email{dremenik@dim.uchile.cl}

\maketitle

\begin{abstract} We review the Airy processes; their formulation and how they are
  conjectured to govern the large time, large distance spatial fluctuations of one
  dimensional random growth models.  We also describe formulas which express the
  probabilities that they lie below a given curve as Fredholm determinants of certain
  boundary value operators, and the several applications of these formulas to variational
  problems involving Airy processes that arise in physical problems, as well as to their
  local behaviour.
\end{abstract}

\tableofcontents

\section{Introduction}
\label{sec:intro}

\subsection{Airy processes and the KPZ universality class}
\label{sec:airykpz}

The \emph{Airy processes} are a collection of stochastic processes which are expected to
govern the long time, large scale, spatial fluctuations of random growth models in the one dimensional
\emph{Kardar-Parisi-Zhang (KPZ) universality class} for wide classes of initial
data. Although there is no precise definition of the KPZ class, it can be identified at
the roughest level by the unusual $t^{1/3}$ scale of fluctuations. It is expected to
contain a large class of random growth processes, as well as randomly stirred one
dimensional fluids, polymer chains directed in one dimension and fluctuating transversally
in the other due to a random potential (with applications to domain interfaces in
disordered crystals), driven lattice gas models, reaction-diffusion models in
two-dimensional random media (including biological models such as bacterial colonies),
randomly forced Hamilton-Jacobi equations, etc. The model giving its name to the
universality class is the \emph{KPZ equation}, which was introduced by \citet{KPZ} as a
model of randomly growing interfaces, and is given by
\begin{equation}
  \label{eq:kpz}
  \p_th=-\tfrac{1}{2}\big(\p_xh\big)^2+\tfrac{1}{2}\p_x^2h+\xi,
\end{equation}
where $\xi(t,x)$ is Gaussian space-time white noise,
$\ee\big(\xi(t,x)\xi(s,y)\big)=\delta_{s=t}\delta_{x=y}$.

A combination of non-rigorous methods (renormalization, mode-coupling, replicas) and
mathematical breakthroughs on a few special models has led to very precise predictions of
universal scaling exponents and exact statistical distributions describing the long time
properties. These predictions have been repeatedly confirmed through Monte-Carlo
simulation as well as experiments; in particular, recent spectacular experiments on
turbulent liquid crystals by Takeuchi and Sano \cite{takeuchiSano1,takeuchiSano2} have
been able to even confirm some of the predicted fluctuation statistics in a physical system.

The conjectural picture that has developed is that the universality class is divided into
subuniversality classes which depend on the class of initial data (or boundary
conditions), but not on other details of the particular models. There are three classes of
initial data which stand out because of their self-similarity properties: Dirac
$\delta_0$, corresponding to curved, or droplet type initial data; $0$, corresponding to
growth off a flat substrate; and $e^{B(x)}$ where $B(x)$ is a two sided Brownian motion,
corresponding to growth in equilibrium. As we will see later, each of these three classes
correspond to concrete initial (or boundary) conditions for the discrete models in the KPZ
class. In addition to these three basic initial data, there are three non-homogeneous
subuniversality classes corresponding roughly to starting with one of the basic three on
one side of the origin, and another on the other side. For one specific discrete model
(last passage percolation or, equivalently, the totally asymmetric exclusion process) the
asymptotic spatial fluctuations have been computed exactly for these six basic classes of
initial data, and are given by the Airy processes: the three basic Airy processes,
Airy$_2$, Airy$_1$ and Airy$_{\rm stat}$, and the crossover Airy processes Airy$_{2\to
  1}$, Airy$_{2\to {\rm BM}}$ and Airy$_{1\to {\rm BM}}$. Although these processes have
been proved to arise as the limiting spatial fluctuations only for one model (and actually
several others in the case of the Airy$_2$ process), as a consequence of the universality
conjecture for the KPZ class it is expected that the same should hold for the other models
in the class.

The purpose of this review is two-fold. First, we will explain in detail in the
 introduction the conjectural picture that we have just sketched from two different
points of view: last passage percolation (or, more generally, directed random polymers)
and the KPZ equation (or, more precisely, the stochastic heat equation). Along the way we
will survey known results for these models.

Our second purpose is to survey a collection of results for the Airy processes which
express the probability that they lie below a given curve as Fredholm determinants of certain
boundary value operators. These expressions have turned out to be very useful in obtaining
some exact distributions through certain variational formulas, and in addition have
allowed one to study some local properties of these processes. This will be the subject of
Sections \ref{sec:fred}-\ref{sec:appl},

\vspace{6pt}
\paragraph{\bf Acknowledgments} JQ was supported by the Natural Science and Engineering
Research Council of Canada. DR was supported by Fondecyt Grant 1120309 and Conicyt
Basal-CMM. The authors thank an anonymous referee for many useful suggestions.

\subsection{Directed random polymers and last passage percolation}
\label{sec:polymersLPP}

\subsubsection{Polymers}Consider the following model of a \emph{directed polymer in a
  random environment}. A \emph{polymer path} is an up-right path $\pi=(\pi_0,\pi_1,\dots)$
in $(\zz_+)^2$ started at the origin, that is, $\pi_0=(0,0)$ and
$\pi_k-\pi_{k-1}\in\{(1,0),(0,1)\}$ (see Figure \ref{fig:lpp}). On $(\zz_+)^2$ we place a
collection of independent random weights $\big\{\omega_{i,j}\big\}_{i,j>0}$. The
\emph{energy} of a polymer path segment $\pi$ of length $N$  is
\[H_N(\pi)=-\sum_{k=1}^N\omega_{\pi_k}.\]
We define the \emph{weight} of such a polymer path segment as
\begin{equation}\label{energy}W_N(\pi)=e^{-\beta H_N(\pi)}=e^{\beta\sum_{k=1}^N\omega_{\pi_k}}\end{equation} for some fixed $\beta>0$ which is known
as the \emph{inverse temperature} of the model.  Let $\Pi_{M,N}$ denote the set of
up-right paths going from the origin to $(M,N)\in(\zz_+)^2$. If we restrict our attention
to such paths then we talk about a \emph{point-to-point polymer}, defined through the
following path measure on $\Pi_{M,N}$:
\begin{equation}
  Q^{\rm point}_{M,N}(\pi)=\frac{1}{Z^{\rm point}(M,N)}W_{M+N}(\pi)\label{eq:Qpoint}
\end{equation}
The normalizing constant \[Z^{\rm point}(M,N)=\sum_{\pi\in\Pi_{M,N}}W_{M+N}(\pi)\] is known
as the \emph{point-to-point partition function}. Similarly, if we consider all possible
paths of length $2N$ then we talk about a \emph{point-to-line polymer}, defined through the
following path measure on $\bigcup_{k=-N,\dots,N}\Pi_{N+k,N-k}$ (that is, all paths of length $2N$):
\begin{equation}
  Q^{\rm line}_N(\pi)=\frac{1}{Z^{\rm line}(N)}W_{2N}(\pi),\label{eq:Qline}
\end{equation}
with the \emph{point-to-line partition function} \[Z^{\rm line}(N)=\sum_{k=-N}^NZ^{\rm
  point}(N+k,N-k).\]

A main quantity of interest in each case is the \emph{free energy}, defined as the
logarithm of the partition function. In the point-to-line case, another important quantity
of interest is the position of the endpoint of the randomly chosen path, which we will
denote by $\kappa_N$. It is widely believed that these quantities should satisfy the scalings
\begin{gather}
  \log(Z^{\rm point}(N,N))\sim a_2N+b_2N^\chi\zeta_2,\label{eq:ZNGUE}\\
  \log(Z^{\rm line}(N))\sim
  a_1N+b_1N^\chi\zeta_1,\label{eq:ZNGOE}\\
  \kappa_N\sim N^\xi\ct\label{eq:ZNPE}
\end{gather}
as $N\to\infty$, where the constants $a_1,a_2$ and $b_1,b_2$ may depend on the distribution of the
$\omega_{i,j}$ and $\beta$, but $\zeta_1,\zeta_2$ and $\ct$ should be universal up to some
fairly generic assumptions on the $\omega_{i,j}$'s, while the fluctuation exponent
\begin{equation}
\chi = 1/3
\end{equation}
and wandering exponent
\begin{equation}
\xi = 2/3.
\end{equation}
Here, and in the rest of this article, whenever we write a relation like
\[Z_N\sim aN+bN^{\kappa}\zeta\] as $N\to\infty$, what we mean is that
\[\lim_{N\to\infty}\pp\!\left(\frac{Z_N-aN}{bN^{\kappa}}\leq m\right)=\pp(\zeta\leq m).\]
One can also have higher, $d+1$ dimensional versions of the model, with the paths directed in
one dimension, and wandering in the other $d$.  In all dimensions the
scaling exponents $\chi$ and $\xi$ are
conjectured to satisfy the \emph{KPZ scaling relation}
\begin{equation}
  \chi=2\xi-1,\label{eq:KPZexponents}
\end{equation}
while the universality of the limiting distributions is unclear except in $d=1$.
For recent progress on  \eqref{eq:KPZexponents}, see 
\cite{chatter,auffDamron1, auffDamron2}.

Although there are few results available in the general case described above, the
zero-temperature limit $\beta\to\infty$, known as \emph{last passage percolation}, is very
well understood, at least for some specific choices of the environment variables
$\omega_{i,j}$. Before introducing this model, we will briefly introduce the Tracy-Widom
distributions from random matrix theory, which will, somewhat surprisingly, play an
important role in the sequel.

\subsubsection{Tracy-Widom distributions}\label{sec:tw}

We will restrict our attention to the distributions arising from the Gaussian Unitary
Ensemble (GUE) and the Gaussian Orthogonal Ensemble (GOE), although these are by no means
the only distributions coming from random matrix theory which appear in the study of
models in the KPZ universality class. The reader can consult \cite{mehta,andGuioZeit} for
good expositions on random matrix theory.

We start with the unitary case. Let $\mathcal{N}(a,b)$ denote a Gaussian random variable
with mean $a$ and variance $b$. An $N\times N$ \emph{GUE matrix} is an (complex-valued)
Hermitian matrix $A$ such that
$A_{i,j}=\mathcal{N}(0,N/\sqrt{2})+\I\mathcal{N}(0,N/\sqrt{2})$ for $i>j$ and
$A_{i,i}=\mathcal{N}(0,N)$.
Here we assume that all the Gaussian variables appearing
in the different entries are independent (subject to the Hermitian condition). The variance
normalization by $N$ was chosen here to make the connection with models in the KPZ
class more transparent. An alternative way to describe the Gaussian Unitary Ensemble is as
the probability measure on the space of $N\times N$ Hermitian matrices $A$ with density
(with respect to the Lebesgue measure on the $N^2$ independent parameters corresponding to
the real entries on the diagonal and the real and imaginary components of the entries
above the diagonal)
\[\frac{1}{Z_N}e^{-\frac1{2N}\tr A^2}\]
for some normalization constant $Z_N$. If $\lambda^N_1,\dotsc,\lambda^N_N$ are the
eigenvalues of such a matrix, then the Wigner semicircle law states that the empirical
eigenvalue density $N^{-1}\sum_{i=1}^N\delta_{\lambda^N_i}$ has approximately a semicircle
distribution on the interval $[-2N,2N]$. The \emph{Tracy-Widom GUE distribution}
\cite{tracyWidom} arises from studying the fluctuations of the eigenvalues of a GUE matrix
at the edge of the spectrum: if we denote by $\lambda^{\rm max}_{\rm GUE}(N)$ the largest
eigenvalue of an $N\times N$ GUE matrix then \cite{tracyWidom}
\[\lambda^{\rm max}_{\rm GUE}(N)\sim2N+N^{1/3}\zeta_2\]
as $N\to\infty$, where $\zeta_2$ has the GUE Tracy-Widom  distribution, which is defined as
follows:
\begin{equation}
  \label{eq:GUE}
  F_{\rm GUE}(s):=\pp(\zeta_2\leq s)=\det(I-P_s\K P_s)_{L^2(\rr)},
\end{equation}
where $K_{\Ai}$ is the \emph{Airy kernel}
\begin{equation}
  K_{\Ai}(x,y)=\int_0^\infty d\lambda\Ai(x+\lambda)\Ai(y+\lambda),\label{eq:airyKernel}
\end{equation}
$\Ai(\cdot)$ is the Airy function, $ P_a$ denotes the projection onto the interval
$(a,\infty)$, and the determinant means the Fredholm determinant on the Hilbert space
$L^2(\rr)$. We will talk at length about the Airy kernel and related operators in later
sections, so for now we will postpone the discussion. Fredholm determinants can be
regarded as the natural generalization of the usual determinant to operators on infinite
dimensional spaces.  We will review their definition and properties in Section
\ref{sec:fred}. Since these determinants will appear often during the rest of this
introduction, the reader who is not familiar with them may want to read Section
\ref{sec:fred} before continuing. 

Before continuing to the $F_{\rm GOE}$ we quickly note that one of the key contributions
of Tracy and Widom \cite{tracyWidom} was to connect \eqref{eq:GUE} to integrable systems.  Let $q(s)$ be the
Hastings-McLeod solution of the Painlev\'{e} II equation
\begin{equation}\label{in:10}
	q''(s)=2q(s)^3+sq(s),
\end{equation}
defined by the additional boundary condition
\begin{equation}\label{in:11}
	q(s) \sim \Ai(s) \quad \textrm{as} \ s \to \infty.
\end{equation}
Then 
\begin{equation}
 F_{\rm GUE}(s) =e^{ -\int_s^\infty dx\hspace{0.05em}(x-s)^2 q^2(x) }.
\end{equation}

The story for the Gaussian Orthogonal Ensemble is similar. An $N\times N$ \emph{GOE
  matrix} is a (real-valued) symmetric matrix $A$ such that
$A_{i,j}=\mathcal{N}(0,N)$ for $i>j$ and $A_{i,i}=\mathcal{N}(0,\sqrt{2}N)$, where
as before we assume that all the Gaussian variables appearing in the different entries are
independent (subject to the symmetry condition). Analogously to the GUE case, the Gaussian
Orthogonal Ensemble can be regarded as the probability measure on the space of $N\times N$
real symmetric matrices $A$ with density
\[\frac{1}{Z_N}e^{-\frac1{4N}\tr A^2}\]
for some normalization constant $Z_N$. As for the GUE, the Wigner semicircle law states
that the empirical eigenvalue density for the GOE has approximately a semicircle
distribution on the interval $[-2N,2N]$. The fluctuations of the spectrum at its edge now
give rise to the \emph{Tracy-Widom GOE distribution}: we denote by $\lambda^{\rm max}_{\rm
  GOE}(N)$ the largest eigenvalue of an $N\times N$ GOE matrix, then \cite{tracyWidom2}
\[\lambda^{\rm max}_{\rm GOE}(N)\sim2N+N^{1/3}\zeta_1\]
as $N\to\infty$, where $\zeta_1$ has the GOE Tracy-Widom  distribution, defined as
\begin{equation}
  \label{eq:GOE}F_{\rm GOE}(s):=\pp(\zeta_1\leq m)=\det(I-P_0B_sP_0)_{L^2(\rr)},
\end{equation}
where $B_s$ is the kernel
\begin{equation}
  B_s(x,y)=\Ai(x+y+s).\label{eq:defB0}
\end{equation}
 This Fredholm determinant formula for
$F_{\rm GOE}$ is essentially due to \cite{sasamoto}, and was proved in
\cite{ferrariSpohn}. The original formula derived by Tracy and Widom is\begin{equation}
  \label{eq:GOE}F_{\rm GOE}(s)=e^{ -\frac12 \int_s^\infty dx\hspace{0.1em}q(x) }  \sqrt{F_{\rm GUE}(s)}
\end{equation}
with $q$ as above.
 
 
\subsubsection{Last passage percolation}\label{sec:lpp}

We come back now to our discussion about directed random polymers, and in particular their
zero-temperature limit. We will restrict the discussion to \emph{geometric last passage
  percolation (LPP)}, where one considers a family $\big\{\omega_{i,j}\}_{i,j>0}$ of
independent geometric random variables with parameter $q$
(i.e. $\pp(\omega_{i,j}=k)=q(1-q)^{k}$ for $k\geq0$). For convenience we also set for now
$\omega_{i,j}=0$ if $i$ or $j$ is 0. As $\beta\to\infty$, the random path measures in
\eqref{eq:Qpoint} and \eqref{eq:Qline} assign an increasingly larger mass to the path
$\pi$ of length $K>0$ which maximizes the weight $W_{K}(\pi)$. In the limit, the path
measures $Q^{\rm point}_{M,N}$ and $Q^{\rm line}_N$ concentrate on the maximizing path,
and the quantities which play the role of the free energy are the \emph{point-to-point
  last passage time},
\[L^{\rm point}(M,N)=\max_{\pi\in\Pi_{M,N}}\sum_{i=0}^{M+N}\omega_{\pi_i}\] and the
\emph{point-to-line last passage time} by
\begin{equation}
L^{\rm line}(N)=\max_{k=-N,\dots,N}L^{\rm point}(N+k,N-k).\label{eq:lpplptline}
\end{equation}
Observe that these last passage times are random, as they depend on the random environment
defined by the $\omega_{i,j}$.

\begin{figure}
  \centering
  \begin{subfigure}[b]{0.35\textwidth}
    \centering 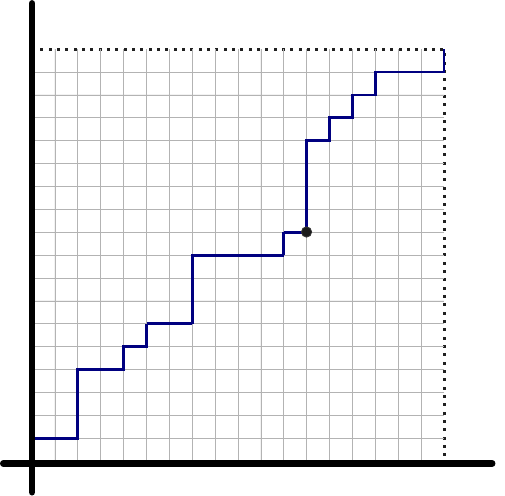
    \subcaption{\label{fig:lpp}}
  \end{subfigure}
  \hskip0.6in
  \begin{subfigure}[b]{0.35\textwidth}
    \centering 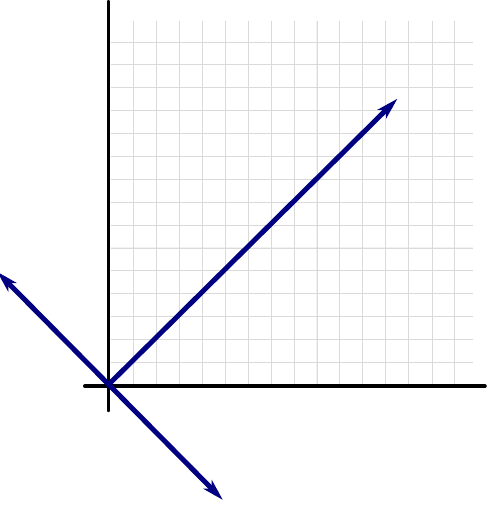
    \subcaption{\label{fig:timespace}}
  \end{subfigure}
  \caption{(a) A polymer/LPP path $\pi$ connecting the origin to $(M,N)$. ~~(b)~Time $s$ and space $u$ axes in LPP.}
\end{figure}

The breakthrough, which in a sense got the whole field started, was the surprising 1999
result by \citet{baikDeiftJohansson} which proved that the asymptotic fluctuations of the
longest increasing subsequence of a random permutation have the Tracy-Widom GUE
distribution. There is an intimate (and simple) connection between this model and LPP
which we will not discuss, instead we will state the companion result  by \citet{johanssonShape} for the point-to-point LPP case:
\begin{equation}
  \label{eq:LPPGUE}
  L^{\rm point}(N,N)\sim c_1N+c_2N^{1/3}\zeta_2,
\end{equation}
where $c_1,c_2$ are some explicit constants which depend only on $q$ and can be found in
\cite{johanssonShape} and $\zeta_2$ has the Tracy-Widom GUE distribution. A similar result
holds for point-to-line LPP. The longest increasing subsequence version goes back to
\citet{baikRains}, while the analogue for LPP which we state here was first proved in
\cite{borFerPrahSasam} (see also \cite{sasamoto,bfp}):
\begin{equation}
  \label{eq:LPPGOE}
  L^{\rm line}(N)\sim c_1'N+c_2'N^{1/3}\zeta_1,
\end{equation}
where $\zeta_1$ now has the GOE Tracy-Widom distribution.

The reason why these exact results (and others we will discuss below) can be obtained for
geometric last passage percolation and other related models is that LPP has an extremely
rich algebraic structure which allows one to write explicit formulas for the distribution of
the last passage times. The algebraic structure arises from regarding the model as a
randomly growing Young tableau, where the cell $(i,j)$ is added at time $L^{\rm
  point}(i,j)$. This shift of perspective relates the problem to the representation theory
of the symmetric group, and in particular to the Robinson-Schensted-Knuth (RSK)
correspondence, which is the main combinatorial tool used in \cite{johanssonShape} to
prove the following remarkable formula:
\[\pp(L^{\rm point}(M,N)\leq s)=\det(I-P_sK^{\rm Meix}_NP_s)_{L^2(\rr)}\]
for $M\leq N$, where the \emph{Meixner kernel} $K^{\rm Meix}_N$ is given by
\[K^{\rm
  Meix}_N(x,y)=\frac{\kappa_N}{\kappa_{N-1}}\frac{p_N(x)p_{N-1}(y)-p_{N-1}(x)p_N(y)}{x-y}\sqrt{w(x)w(y)},\]
$w(x)=\binom{M-N+x}{x}$, and the functions $p_N(x)$ are the normalized Meixner
polynomials, i.e., the normalized family of discrete orthogonal polynomials $p_N(x)$ with respect
to the weight $w(x)$, with $p_N(x)$ of degree $N$ and leading coefficient $\kappa_N$. A
non-trivial asymptotic analysis of this kernel allowed Johansson to deduce that the above
Fredholm determinant converges as $N\to\infty$ to the Fredholm determinant appearing in
the definition \eqref{eq:GUE} of the Tracy-Widom GUE distribution. A more detailed
discussion of these facts is beyond the scope of this review; what the reader should keep
from this discussion is that the exact results which we are discussing depend crucially on
what is usually referred to as \emph{exact solvability} or \emph{integrability}: the availabilty of (extremely
non-trivial) exact formulas for quantities of interest. These formulas arise from the very rich
algebraic structure present in some (but by no means all) models in the KPZ class. For a
recent survey on this subject see \cite{borodinGorinReview}.

As part of the general KPZ universality conjecture, one expects that
\eqref{eq:LPPGUE} and \eqref{eq:LPPGOE} hold not only for LPP, but in general for any
$\beta>0$. In other words, the belief is that in \eqref{eq:ZNGUE} and \eqref{eq:ZNGOE},
the random variables $\zeta_2$ and $\zeta_1$ have respectively the Tracy-Widom GUE and GOE
distributions. There has been only partial progress in proving this conjecture for
point-to-point directed polymers (and virtually none in the point-to-line case), the
difficulty lying in the lack of exact solvability. Versions of this conjecture have been
proved for two related models in the point-to-point case: the continuum random polymer in
\cite{acq} (building on results of \cite{tracyWidomASEP2,tracyWidomASEP1,tracyWidomASEP3})
and the semi-discrete polymer of O'Connell and Yor in \cite{borCor,borodinCorwinFerrari}
(see also \cite{oconnellQTL}). In the setting of discrete directed random polymers,
\cite{cosz} showed that if the weights are chosen so that $-w_{i,j}$ is distributed as the
logarithm of a Gamma random variable with parameter $\theta_i+\hat\theta_j$ (for some
fixed $\theta_i$'s and $\hat\theta_j$'s) then the model is exactly solvable in the sense
explained above. This was later used in \cite{bcrLogGamma} to prove that the asymptotic
fluctuations of the free energy of the point-to-point polymer (at least for low enough
temperature) have the conjectured Tracy-Widom GUE distribution.

\subsubsection{Spatial fluctuations and the Airy processes}\label{sec:spatialfluc}

The Airy processes arise from LPP when we look not only at the fluctuations of the free
energy at a single site, but instead at several sites. To this end, we define the rescaled
point-to-point process $u\mapsto H^{\rm point}_N(u)$ by linearly interpolating the values
given by scaling $L^{\rm point}(M,N)$ through the relation
\begin{equation}
L^{\rm point}(N+u,N-u)=c_1N+c_2N^{1/3}H^{\rm point}_N(c_3N^{-2/3}u)\label{eq:rescaledptp}
\end{equation}
for $u=-N,\dotsc,N$, where the constants $c_i$ have explicit expressions which depend only on
$q$ and can be found in \cite{johansson}. Observe that this corresponds to looking at the
free energy at a line of slope $-1$ passing through $(N,N)$. The limiting behavior of
$H^{\rm point}_N$ is described by the Airy$_2$ process $\aip$ (minus a parabola, see
Theorem \ref{thm:Airy2LPP}). This  process was
introduced by \citet{prahoferSpohn}, and is defined through its finite-dimensional
distributions, which are given by a Fredholm determinant formula: given
$x_0,\dots,x_n\in\mathbb{R}$ and $u_1<\dots<u_n$ in $\mathbb{R}$,
\begin{equation}\label{eq:detform}
  \mathbb{P}\!\left(\aip(u_1)\le x_1,\dots,\aip(u_n)\le x_n\right) =
  \det(I-\mathrm{f}^{1/2}\K^{\mathrm{ext}}\mathrm{f}^{1/2})_{L^2(\{u_1,\dots,u_n\}\times\mathbb{R})},
\end{equation}
where we have counting measure on $\{u_1,\dots,u_n\}$ and Lebesgue measure on
$\mathbb{R}$, $\mathrm f$ is defined on $\{u_1,\dots,u_n\}\times\mathbb{R}$ by
\begin{equation}
  \mathrm{f}(u_j,x)=\uno{x\in(x_j,\infty)}, \label{eq:deff}
\end{equation}
and the {\it extended Airy kernel} \cite{prahoferSpohn,FNH,macedo} is defined by
\begin{equation}
  \K^\mathrm{ext}(u,\xi;u',\xi')=
  \begin{cases}
    \int_0^\infty d\lambda\,e^{-\lambda(u-u')}\Ai(\xi+\lambda)\Ai(\xi'+\lambda), &\text{if $u\ge u'$}\\
    -\int_{-\infty}^0 d\lambda\,e^{-\lambda(u-u')}\Ai(\xi+\lambda)\Ai(\xi'+\lambda),
    &\text{if $u<u'$}.
  \end{cases}\label{eq:extAiry}
\end{equation}
Although it is not obvious from the definition, the Airy$_2$ process is stationary (this
will become clear in Section \ref{sec:bdVlK}), and as should be expected from
\eqref{eq:LPPGUE}, $\pp(\aip(u)\leq m)=F_{\rm GUE}(m)$ for all $u$.  There is a close
connection, which we will explain in Section \ref{sec:extendedDet}, between the Airy
kernel $\K$ appearing in the definition \eqref{eq:GUE} of the Tracy-Widom GUE distribution
and the extended kernel $\K^{\rm ext}$.

The precise result linking the point-to-point LPP spatial fluctuations to the Airy$_2$
process is due to \citeauthor{johansson} (see also \cite{prahoferSpohn}):

\begin{thm}[\citet{johansson}]\label{thm:Airy2LPP}
  There is a continuous version of $\aip$, and
  \[H^{\rm point}_N(u)\xrightarrow[N\to\infty]{}\aip(u)-u^2\] in distribution in the
  topology of uniform convergence of continuous functions on compact sets.
\end{thm}

In the LPP picture, the ``time'' variable (which we will denote by $s$) flows in the
$(1,1)$ direction of the plane, while ``space'' (which, as above, we will denote by $u$)
corresponds to the direction $(1,-1)$ (see Figure \ref{fig:timespace}). In this sense,
Theorem \ref{thm:Airy2LPP} describes the spatial fluctuations of the point-to-point last
passage times as time $s\to\infty$.

\begin{figure}
  \centering
  \begin{subfigure}[b]{0.4\textwidth}
    \centering 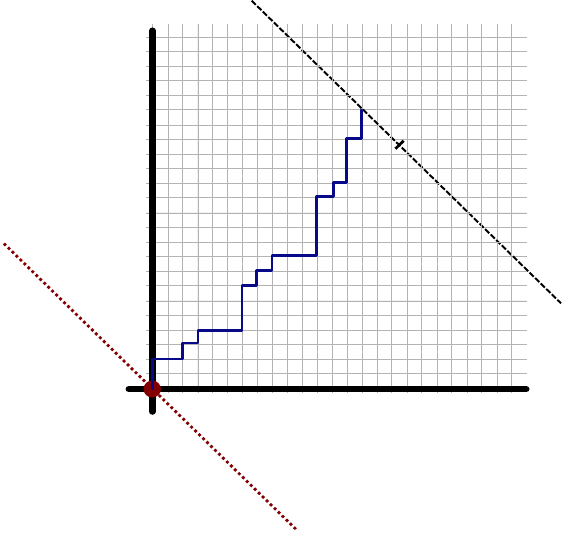
    \subcaption{Narrow wedge\label{fig:wedge}}
  \end{subfigure}
  \hskip0.6in
  \begin{subfigure}[b]{0.4\textwidth}
    \centering 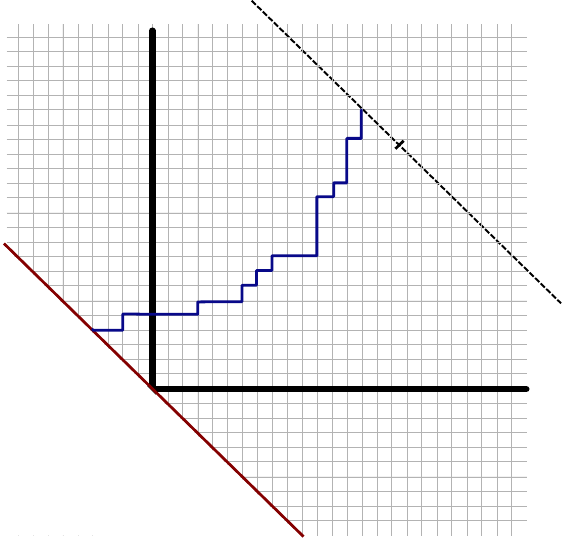
    \subcaption{Flat\label{fig:flat}}
  \end{subfigure}
  \hskip0.3in
  \begin{subfigure}[b]{0.4\textwidth}
    \centering 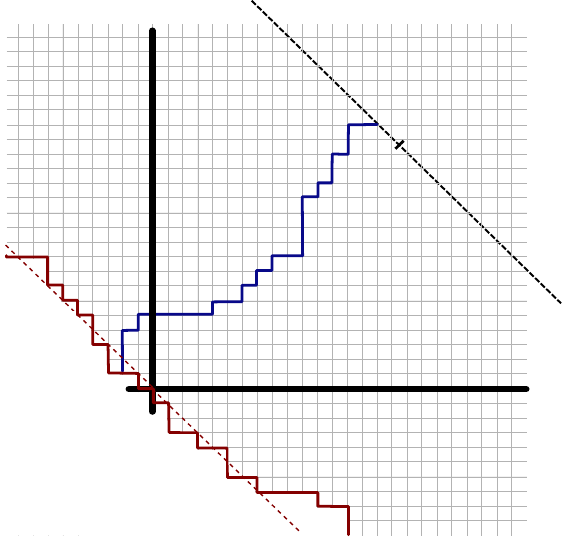
    \subcaption{Stationary\label{fig:stat}}
  \end{subfigure}
  \caption{Schematic representation of the LPP models with asymptotic spatial fluctuations
    given by: (a) Airy$_2$; (b) Airy$_1$; (c) Airy$_{\rm stat}$.}
\end{figure}

One can think of extending the LPP model to paths starting at $s=0$ with any space
coordinate, i.e., paths which start at any point of the form $(k,-k)$, $k\in\zz$. To
recover point-to-point LPP one simply sets $\omega_{i,j}=0$ whenever $i\leq0$ or $j\leq0$,
which is easily seen to be equivalent (from the point of view of last passage times) to
forcing our paths to start at the origin. In this sense, point-to-point LPP and the
Airy$_2$ process correspond to the $\delta_0$ (also known as \emph{delta}, \emph{narrow
  wedge} or \emph{curved}) initial data (see Figure \ref{fig:wedge}). Note that in this
case we only assign positive weights to sites such that $s>|u|$. To recover the flat,
stationary and mixed initial data which we introduced earlier, we need to assign weights
to sites such that $s\leq|u|$.

\begin{rem}\label{rem:diff}
  The results for the flat, stationary and mixed initial data have been proved in settings
  which differ slightly from the one introduced here. To avoid additional notation and
  complications, we will state the results on the present setting. We refer the reader to
  the corresponding references for more details on the differences. In the case of
  multipoint results, one can translate between the various settings by using the slow
  decorrelation result proved of \cite{ferrari-slow,corwinFerrariPeche1}, as done in
  \cite{baikFerrPeche,corwinFerrariPeche}.
\end{rem}

We start with the \emph{flat} initial data. It corresponds to extending the weights
$\omega_{i,j}$ to be independent geometric random variables with parameter $q$ whenever $i+j>0$
and setting $\omega_{i,j}=0$ otherwise (see Figure \ref{fig:flat}). This corresponds
to letting our paths start at any site in the line $s=0$ but not attach any additional
weights along that line, which explains the name flat. The corresponding point-to-line
rescaled process may be defined as follows: first we extend the definition of last passage
times to accomodate the flat initial data,
\[L^{\rm point}_{\rm
  flat}(M,N)=\max_{i\in\zz}\max_{\pi\in\Pi_{(i,-i)\to(M,N)}}\sum_{j=0}^{2i+M+N}\omega_{\pi_j}\]
with self-explanatory notation, and then we define the rescaled process $u\mapsto H^{\rm
  line}_N(u)$ by linearly interpolating the values given the relation
\[L^{\rm point}_{\rm flat}(N+u,N-u)=c_1N+c_2N^{1/3}H^{\rm line}_N(c_3N^{-2/3}u)\] for
$u=-N,\dotsc,N$.  The flat initial data gives rise to the Airy$_1$ process $\aipo$, which
was introduced by \citet{sasamoto}, and is defined through its finite-dimensional
distributions,
\begin{equation}\label{eq:detAiry1}
  \pp\!\left(\aipo(u_1)\le \xi_1,\dots,\aipo(u_n)\le \xi_n\right) = 
  \det\!\big(I-\mathrm{f}K^{\mathrm{ext}}_1\mathrm{f}\big)_{L^2(\{u_1,\dots,u_n\}\times\mathbb{R})},
\end{equation}
with ${\rm f}$ as in \eqref{eq:deff} and
\begin{multline}\label{eq:fExtAiry1}
  K^{\rm ext}_1(u,\xi;u',\xi')=-\frac{1}{\sqrt{4\pi
      (u'-u)}}\exp\!\left(-\frac{(\xi'-\xi)^2}{4 (u'-u)}\right)\uno{u'>u}\\
  +\Ai(\xi+\xi'+(u'-u)^2) \exp\!\left((u'-u)(\xi+\xi')+\tfrac23(u'-u)^3\right).
\end{multline}
The Airy$_1$ process is stationary, and as should be
expected from \eqref{eq:LPPGOE}, its marginals are given by the Tracy-Widom GOE
distribution: $\pp(\aipo(u)\leq m)=F_{\rm GOE}(2m)$ for all $u$.

\begin{thm}[\cite{borFerPrahSasam,bfp,borodFerSas}]
  \[H^{\rm line}_N(u)\xrightarrow[N\to\infty]{}2^{1/3}\aipo(2^{-2/3}u)\] in the sense of
  convergence of finite-dimensional distributions (on a slightly different setting than
  the one presented here, see Remark \ref{rem:diff}).
\end{thm}

The powers of $2^{1/3}$ in the above limit should be regarded as an arbitrary
normalization (in fact, one could have defined the Airy$_1$ process as this scaled
version of it). The appearence of these factors has to do with the fact that the natural scaling
in the definition of these quantities differs between random matrix models and models such
as LPP or directed polymers (as, for instance, in \eqref{eq:Airy2GOE} below). See Section
\gref{sec:indirect} for a related discussion.

The \emph{stationary} initial data is slightly more cumbersome to introduce. The name
stationary comes from the fact that for the closely related totally asymmetric exclusion
process (TASEP), this initial condition corresponds to starting with particles placed
according to a product Bernoulli measure with parameter $1/2$, which is stationary for the
process. Translated to LPP, this initial condition corresponds to the following. Let
$(S_n)_{n\in\zz}$ be the path of a double-sided simple random walk on $\zz$ with $S_0=0$,
which we assume to be independent of the weights $\omega_{i,j}$. We rotate this random
walk path by an angle of $-\pi/4$ and then put it along the $s=0$ line by defining the
(random) discrete curve
$\gamma_0=\big\{\big(\frac12(S(i)+i),\frac12(S(i)-i)\big),\,i\in\zz\big\}$. We then extend
the weights $\omega_{i,j}$ to be independent geometric random variables with parameter $q$
whenever $(i,j)$ lies above $\gamma_0$ and $\omega_{i,j}=0$ otherwise (see Figure
\ref{fig:stat}). The corresponding stationary rescaled process $H^{\rm stat}_N(u)$ can be
defined analogously to the previous cases, by maximizing over paths starting at $\gamma_0$
and going to the anti-diagonal line passing through $(N,N)$. It gives rise to the
Airy$_{\rm stat}$ process $\mathcal{A}_{\rm stat}$. Its definition is also given in terms
of finite-dimensional distributions involving Fredholm determinants, but the formulas are
a lot more cumbersome. We will not need the exact formulas, so we refer the reader to
\cite{baikFerrPeche} for the details. Despite its name, $\mathcal{A}_{\rm stat}$ is not
stationary as a process. In fact, due to the connection with stationary TASEP,
$\mathcal{A}_{\rm stat}$ is just a standard double-sided Brownian motion, but with a
non-trivial random height shift at the origin given by the Baik-Rains distribution, see
\cite{baikRainsF0}. The convergence result in this case is the following:

\begin{thm}[\cite{baikFerrPeche}]
  \[H^{\rm stat}_N(u)\xrightarrow[N\to\infty]{}\mathcal{A}_{\rm stat}(u)\] in the sense of
  convergence of finite-dimensional distributions (on a slightly different setting than
  the one presented here, see Remark \ref{rem:diff}).
\end{thm}

The mixed initial conditions can be obtained by placing one condition on each half of the
line $u=0$. We will explain how this is done in the case of the \emph{half-flat}, or
\emph{wedge$\to$flat} initial data, and leave the examples leading to $\mathcal{A}_{2\to
  \rm BM}$ and $\mathcal{A}_{1\to \rm BM}$ to the interested reader (see
\cite{bfsTwoSpeed,corwinFerrariPeche}). To obtain the Airy$_{2\to1}$ process we extend the
weights $\omega_{i,j}$ to be independent geometric random variables with parameter $q$ whenever
$i,j>0$, or $i+j>0$ with $i<0$, setting $\omega_{i,j}=0$ for all other sites. The half-flat
rescaled process $H^{\rm half-line}_N(u)$ is obtained as in the previous cases, and gives
rise to the Airy$_{2\to1}$ process $\Bt$. It was introduced by \citet{bfs}, and is given
by
\begin{equation}\label{eqTransAiryProcess}
  \pp\!\left(\Bt(u_1)\le\xi_1,\dots,\Bt(u_m)\le\xi_m\right)=
  \det\!\big(I-{\rm f} K^{\rm ext}_{2\to1}{\rm f}\big)_{L^2(\{u_1,\dots,u_m\}\times\mathbb{R})},
\end{equation}
with ${\rm f}$ as in \eqref{eq:deff} and
\begin{multline}\label{eqKCompleteInfinity}
  K^{\rm ext}_{2\to1}(u,\xi;u',\xi')=-\frac{1}{\sqrt{4\pi(\xi'-\xi)}}\exp\!\left(-\frac{(\tilde\xi'-\tilde\xi)^2}{4(u'-u)}\right)\uno{u'>u}\\
  +\frac{1}{(2\pi \I)^2} \int_{\gamma_+}dw \int_{\gamma_-}dz\, \frac{e^{w^3/3+u' w^2-\tilde\xi'
      w}}{e^{z^3/3+u z^2-\tilde\xi z}} \frac{2w}{(z-w)(z+w)},
\end{multline}
where $\tilde\xi=\xi-u^2\uno{u\leq0}$, $\tilde\xi'=\xi'-(u')^2\uno{u'\leq 0}$ and the paths
$\gamma_+,\gamma_-$ satisfy $-\gamma_+\subseteq\gamma_-$ with
\mbox{$\gamma_+:e^{\I\phi_+}\infty\to e^{-\I\phi_+}\infty$},
\mbox{$\gamma_-:e^{-\I\phi_-}\infty\to e^{\I\phi_-}\infty$} for some $\phi_+\in
(\pi/3,\pi/2)$, $\phi_-\in (\pi/2,\pi-\phi_+)$. As could be expected from the above
description, the Airy$_{2\to 1}$ process crosses over between the Airy$_2$ and the
Airy$_1$ processes in the sense that $\Bt(u+v)$ converges to
$2^{1/3}\aipo(2^{-2/3}u)$ as $v\to\infty$ and $\aip(u)$ when $v\to-\infty$. The
convergence result is the following:

\begin{thm}[\cite{bfs}]
  \[H^{\rm half-line}_N(u)-u^2\uno{u\leq0}\xrightarrow[N\to\infty]{}\mathcal{A}_{2\to1}(u)\] in
  the sense of convergence of finite-dimensional distributions (on a slightly different setting than
  the one presented here, see Remark \ref{rem:diff}).
\end{thm}

Some of these results have been extended to the case where the points at which one computes the
corresponding finite-dimensional distributions do not all lie in the same anti-diagonal
line, but instead fall in certain \emph{space-like curves} lying close
enough to such a line, see e.g. \cite{bp-push,borodFerSas} for more details and
\cite{corwinFerrariPeche} for further extensions.

From the definitions it is clear that the basic three Airy processes $\mathcal{A}_{2}$, $\mathcal{A}_{1}$,
and  $\mathcal{A}_{\rm stat}$ are invariant  under $\mathcal{A}(u) \mapsto \mathcal{A}(-u)$, but the 
mixed cases are not.

Since all initial data are superpositions of Dirac masses, there is a sense in which the
Airy$_2$ process is the most basic. For example, using the fact that point-to-line last
passage times are computed simply as the maximum of point-to-point last passage times,
\citet{johansson} obtained the following celebrated formula as a corollary of
\eqref{eq:LPPGUE} and Theorem \ref{thm:Airy2LPP}:
\begin{equation}
  \pp\!\left(\sup_{x\in\rr}\big\{\aip(x)-x^2\big\}\leq m\right)=F_{\rm GOE}(4^{1/3}m).\label{eq:Airy2GOE}
\end{equation}
A direct proof of this formula was later provided in \cite{cqr}. The argument used in this
second proof starts with a different expression for the finite-dimensional distributions
of $\aip$ in terms of the Fredholm determinant of a certain boundary value operator. This
type of formula, and their extensions to continuum statistics, are the starting point of
most of the results we will survey in Sections \ref{sec:contStat} and \ref{sec:appl}. For
example, as described in Section \ref{sec:endpoint}, they allow to compute the asymptotic
distribution of $\kappa_N$, the position of the endpoint in the maximizing path in point-to-line LPP.

Extrapolating from \eqref{eq:Airy2GOE} leads to a conjecture that the one-point
marginals of the other Airy processes should be obtained through certain variational
problems involving the Airy$_2$ process. To state the precise conjectures we turn to the
stochastic heat equation, whose logarithm is the solution of the KPZ equation. The
advantage of this model over LPP and other discrete models is that it is linear in the
initial data, and hence the heuristics are more easily stated in that context.  The disadvantage
is that most of the argument relies on conjectures based on universality.

\subsection{The continuum random polymer and the stochastic heat equation}
\label{sec:she}

We now consider the continuum version of the finite temperature discrete random
polymers \eqref{eq:Qpoint} and \eqref{eq:Qline}.  
The (point-to-point) continuum random polymer is a random probability measure $P^{\beta,\xi}_{T,x}$ on continuous functions $x(t)$
on $[0,T]$ with $x(0)=0$ and $x(T)=x$  with formal weight 
\begin{equation}
e^{-\beta\int_0^Tdt\,\xi(t,x(t)) - \tfrac12\int_0^Tdt\hspace{0.05em}|\dot{x}(t)|^2 }
\end{equation}
given to the path $x(\cdot)$, where $\xi(t,x)$, $t\ge 0$, $x\in \mathbb{R}$ is space-time
white noise, i.e. the distribution-valued Gaussian variable such that for smooth functions
$\varphi$ of compact support in $\mathbb{R_+}\times \mathbb{R}$, $\langle\varphi,
\xi\rangle := \int_{\mathbb{R}_+\times \mathbb{R}}dt\,dx\, \varphi(t,x) \xi(t,x)$ are mean
zero Gaussian random variables with covariance structure $E[ \langle\varphi_1,
\xi\rangle\langle\varphi_1, \xi\rangle] = \langle\varphi_1,\varphi_2\rangle$.  One can
also think of the continuum random polymer as having a density
 \begin{equation}
e^{-\beta\int_0^T dt\,\xi(t,x(t))}
\end{equation}
with respect to the Brownian bridge. Neither prescription makes mathematical sense, but
the second one does if one smooths out the white noise $\xi(t,x)$ in space. Removing the
smoothing, one find that there is indeed a limiting measure supported on continuous
functions $C[0,T]$ which we call $P^{\beta,\xi}_{T,x}$.  In fact, it is a Markov process,
and one can define it directly as follows.  Let $z(s,x,t,y)$ denote the solution of the
stochastic heat equation after time $s\ge 0$ starting with a delta function at $x$,
\begin{equation} 
 \partial_t z =\tfrac12\partial_y^2z - \beta\xi z, \qquad t>s, ~y\in\mathbb{R},\label{shes} \qquad z(s,x,s,y) =\delta_x(y) . 
\end{equation}
It is important that they are all using the same noise $\xi$.  Note that the stochastic
heat equation is well-posed \cite{walsh}.  The solutions look locally like exponential
Brownian motion in space.  They are H\"older $\frac12-\delta$ for any $\delta>0$ in $x$
and $\frac14-\delta$ for any $\delta>0$ in $t$. In fact, exponential Brownian motion
$e^{B(x)}$ is invariant up to multiplicative constants, i.e. if one starts \eqref{shes}
with $e^{B(x)}$ where $B(x)$ is a two-sided Brownian motion, then there is a (random)
$C(t)$ so that $C(t) z(t,x)$ is an exponential of another two-sided Brownian motion
\cite{berGiaco}.  $P^{\beta,\xi}_{T,x}$ is then defined to be the probability measure on
continuous functions $x(t)$ on $[0,T]$ with $x(0)=0$ and $x(T)=x$ and finite dimensional
distributions
\begin{multline}
P^{\beta,\xi}_{T,x}( x(t_1)\in dx_1, \ldots, x(t_n)\in dx_n) \\ = \frac{ z(0,0,t_1,x_1)z(t_1,x_1,t_2,x_2)\cdots z(t_{n-1},x_{n-1},t_n,x_n)z(t_n,x_n,T,x)}{z(0,0,T,x)} dx_1\cdots dx_n
\end{multline}
for $0<t_1<t_2<\cdots<t_n <T$.  One can check these are a.s. a consistent family of finite
dimensional distributions. This holds basically because of the Chapman-Kolmogorov equation
\begin{equation}\label{CK}
\int_{-\infty}^\infty du\, z(s,x,\tau,u)z(\tau,u, t,y) = z(s,x,t,y)
\end{equation}
for $s<\tau\leq t$, which is a consequence of the linearity of the stochastic heat equation.  

Note that the construction is for each $T>0$ fixed.  Unlike the usual case of diffusions,
the measures are very inconsistent for varying $T$.  One should imagine that the polymer
paths are peeking into the future to see the best route, so the measure depends
considerably on all the noise in the time interval $[0,T]$.  We can also define the joint
measure $\mathbb{P}^{\beta}_{T,x}= P^{\beta,\xi}_{T,x}\otimes Q(\xi)$ where $Q$ is the
distribution of the $\xi$, i.e. the probability measure of the white noise.

\begin{thm}[\cite{akq1}]\mbox{}
  \begin{enumerate}[label=(\roman*)]
  \item  The measures $P^{\beta,\xi}_{T,x}$ and $\mathbb{P}^{\beta}_{T,x}$ are
    well-defined (the former, $Q$-almost surely).
  \item $P^{\beta,\xi}_{T,x}$ is a Markov
    process supported on H\"older continuous functions of exponent $\frac12-\delta$ for any
    $\delta>0$, for $Q$-almost every $\xi$.
  \item Let $t_k^n = \frac{k}{2^n}$. Then with
    $\mathbb{P}^{\beta}_{T,x}$ probability one, we have that for all $0 \leq t \leq 1$
    \[\sum_{k=1}^{\lfloor 2^n t \rfloor} \left( x(t_k^n) - x(t_{k-1}^n) \right)^2 \xrightarrow[n\to\infty]{} t,\]
    i.e.  the quadratic variation exists, and coincides with the one obtained for
    $\mathbb{P}^{0}_{T,x}$ (the Brownian bridge measure).
    \item $P^{\beta,\xi}_{T,x}$ is {\it
      singular} with respect to $\pp^{0}_{T,x}$ (the Brownian bridge measure) for $Q$-almost every
    $\xi$.
  \end{enumerate}
\end{thm}

So the continuum random polymer looks locally like, but is singular with respect to,
Brownian motion.  One can also define the point-to-line continuum random polymer
$\mathbb{P}^{\beta}_{T}$, in the same way as in the discrete case.  For large $T$, one
expects ${\rm Var}_{\mathbb{P}^{\beta}_{T}}(x(T))\sim T^{4/3}$ in the point-to-line case
or ${\rm Var}_{\mathbb{P}^{\beta}_{T,0}}(x(T/2))\sim T^{4/3}$ in the point-to-point case.
Here the variance is over the random background as well as $P^{\beta,\xi}_{T,x}$.  The
conditional variance given $\xi$ should be much smaller.

If $z(t,x)$ is the solution of (\ref{shes}) then $h(t,x)= -\beta^{-1}\log z(t,x)$ can be
thought of as either the (renormalized) free energy of the point-to-point continuum random
polymer, or the {\it Hopf-Cole solution of the Kardar-Parisi-Zhang equation},
\begin{equation}\label{KPZ}
\partial_t h = -\tfrac{\beta}2(\partial_x h)^2 +\tfrac12 \partial_x^2 h + \xi,
\end{equation}
for random interface growth. Since $\log z(t,x)$ looks locally like Brownian motion,
\eqref{KPZ} is not well-posed (see \cite{hairer} for recent progress on this question.)
If $\xi$ were smooth, then the Hopf-Cole transformation takes \eqref{shes} to \eqref{KPZ}.
For white noise $\xi$, we take $h(t,x)=-\beta^{-1} \log z(t,x)$ with $z(t,x)$ a solution
of \eqref{shes} to be the {\it definition} of the solution of \eqref{KPZ}.  It is known \cite{berGiaco}
that these are the solutions one obtains if one smooths the noise, solves the equation,
and takes a limit as the smoothing is removed (and after subtraction of a diverging
constant).  They are also the solutions obtained as the limit of discrete models like
asymmetric exclusion in the weakly asymmetric limit \cite{berGiaco}, or directed polymers in the {\it
  intermediate disorder limit} \cite{akq2}.

To understand the intermediate disorder limit we consider how the KPZ equation \eqref{KPZ}
rescales.  Let
\begin{equation}\label{scal}
h_\e (t,x) = \e^a h(\e^{-z} t, \e^{-1} x)\end{equation}
Recall the white noise has the distributional scale invariance
\begin{equation}
\xi(t,x) \stackrel{\rm dist}{=} \e^{\frac{z+1}{2}} \xi(\e^{z} t,\e^{1} x).
\end{equation}
Hence, setting $\beta =1$ for clarity,
\begin{equation}\label{scaling}
\partial_t h_\e =-\tfrac12 \e^{2-z-a} (\partial_x h_\e)^2 + \tfrac12 \e^{2-z} \partial_x^2 h_\e + \e^{a-\frac12 z+\frac12} \xi.
\end{equation}
Because the paths of $h$ are locally Brownian in $x$ we are forced to take $a=1/2$ to see  non-trivial limiting behaviour.   This forces us to take
\begin{equation}
z=3/2
\end{equation}
The non-trivial limiting behaviour of models in the KPZ universality class are all obtained in this scale.

On the other hand, if we started with KPZ with noise of order $\e^{1/2}$,
\begin{equation}\label{scaling2}
\partial_t h =-\tfrac12 (\partial_x h)^2 + \tfrac12  \partial_x^2 h +  \e^{1/2}\xi,
\end{equation}
then a diffusive scaling, 
\begin{equation}\label{scal}
h_\e (t,x) =  h(\e^{-2} t, \e^{-1} x),\end{equation}
would bring us back to the standard KPZ equation (\ref{KPZ}).   This is the intermediate disorder scaling
in which KPZ and the continuum random polymer can be obtained from
discrete directed polymers.  It tells us that if we set 
\begin{equation}
\beta= \e^{1/2}\tilde\beta
\end{equation}
in \eqref{energy}
then the distribution of the rescaled polymer path
\begin{equation}
x_\e (t) := \e x_{\lfloor \e^{-2} t\rfloor}\qquad 0\le t\le T
\end{equation} will converge to the continuum random polymer, with temperature $c\tilde\beta$
(see \cite{akq2} for details).

\subsection{General conjectural picture for the SHE}
\label{sec:conjSHE}
%
%
%

Define  $A_t$ from the solution of \eqref{shes} by
\begin{equation}\label{eq:zed}
  z(0,y; t,x) = \tfrac{1}{\sqrt{2\pi t}} e^{ -\frac{(x-y)^2}{2t} -\frac{t}{24}+ 2^{-1/3}t^{1/3}A_t(2^{-1/3}t^{-2/3}(x-y))}.
\end{equation}
$A_t(\cdot)$ is called the {\it crossover Airy process}, the key conjecture
being \begin{equation}\label{airyconj}A_t( x) \to \aip(x)\end{equation} This is known in
the sense of one-dimensional distributions (see \cite{acq}, where \eqref{airyconj} is
Conjecture 1.5).  A non-rigorous derivation based on a factorization approximation for the
Bethe eigenfunctions of the $\delta$-Bose gas can be found in \cite{prolhacSpohn}. Note
however that the factorization assumption is almost certainly false.

Now one tries to use the linearity of the stochastic heat equation to solve
for general initial data $z(0,x)=z_0(x)$,
\begin{equation}\label{eq:zed1}
  z(t,x) =\int_{-\infty}^\infty dy\, \tfrac{1}{\sqrt{2\pi t}} e^{ -\frac{(x-y)^2}{2t} -\frac{t}{24}+ 2^{-1/3}t^{1/3}A_t(2^{-1/3}t^{-2/3}(x-y))} z_0(y).
\end{equation}
It is not hard to see that the equality is correct in the sense of one-dimensional
distributions, but not more.  If one wants, for example, joint distributions of $z(t,x_i)$
for more than one $x_i$, then one needs to enhance the crossover Airy process in
\eqref{eq:zed} to a two parameter process $A_t(2^{-1/3}t^{-2/3}x, 2^{-1/3}t^{-2/3}y)$.
The conjectural limit of this is a two parameter process we call the {\it Airy sheet}.
However, we do not even have a full conjecture for its finite dimensional distributions,
though some properties can be described (see \cite{corwinQuastelKPZfxpt}).

Calling $\tilde{x}=2^{-1/3}t^{-2/3} x$ and $\tilde{y}=2^{-1/3}t^{-2/3} y$ and  starting with initial data 
$z_0(x)= \exp\{2^{-1/3}t^{1/3}f(2^{-1/3}t^{-2/3}x)  \}$, we can rewrite
the exponent in \eqref{eq:zed1} as
\[2^{-1/3}t^{1/3}\big[A_t (\tilde{x}-\tilde{y})
-(\tilde{x}-\tilde{y})^2-f(\tilde{y}) \big]-\tfrac1{24}t\] so that for large $t$ the fluctuation field
$2^{1/3}t^{-1/3} \big[\log z(t,x)+\frac1{24}t+\log(\sqrt{2\pi t})\big]$ is well
approximated by
\begin{equation}
  \sup_{\tilde{y}\in\rr}\big\{ \aip (\tilde{x}-\tilde{y})-(\tilde{x}-\tilde{y})^2-\tilde{f}(\tilde{y})\big\}.
\end{equation}
The type of initial data would appear to be quite restrictive, but actually this picks out the appropriate
self-similar classes.  The easiest example is the flat case $f=0$.  We obtain the statement
\begin{equation}
  \aipo(x)\quad{\buildrel {\rm (d)} \over =}\quad\sup_{y\in\rr} \big\{\aip(y-x)-(y-x)^2\big\}
  \label{eq:conj}
\end{equation}
in the sense of one-dimensional distributions.  Since the left hand side is just the GOE Tracy-Widom law
 this is the well known theorem of Johansson \eqref{eq:Airy2GOE} once again.   
 
If one starts with a two sided Brownian motion, then the required self-similarity of this initial data is just
the Brownian scaling and one arrives at\footnote{We thank J. Baik and Z. Liu for pointing
  out the missing $\sqrt{2}$ on the right hand side of this equality in an earlier version
  of this manuscript. See \cite{baikLiu} for more details.}
\begin{equation}
 \mathcal{A}_{\rm stat}(x)\quad{\buildrel {\rm (d)} \over =}\quad\sup_{y\in\rr} \big\{\aip(y-x)-(y-x)^2- \sqrt{2}B(y)\big\}.
  \label{eq:conj39}
\end{equation}

The mixed cases require a tiny bit more care.
Let's explain the heuristics first for the case of the Airy$_{2\to1}$ process.  Starting from the step initial data $z(0,x) =\uno{x>0}$ the
prediction is
\begin{equation}\label{eq:1to2approx}
  -\log z(t,x) \approx \tfrac1{2t}x^2\uno{x<0}+\tfrac{1}{24}t+\log(\sqrt{2\pi t}) - 2^{-1/3}t^{1/3} \Bt(2^{-1/3}t^{-2/3} x).
\end{equation}
On the other hand, by linearity we have for each fixed $x$, in distribution,
\begin{equation}\label{eq:zed2}
  z( t,x) =\int_0^\infty dy\,z(0,y; t,x)= \int_0^\infty dy\,\tfrac{1}{\sqrt{2\pi t}}
  e^{ -\frac{(x-y)^2}{2t} -\frac{t}{24}+ 2^{-1/3}t^{1/3}A_t(2^{-1/3}t^{-2/3} (x-y)) }.
\end{equation}
Comparing with \eqref{eq:1to2approx} we deduce that the processes $\sup_{y\ge 0} \big(
\aip(x-y)-(x-y)^2\big)$ and $\Bt(x)-x^2\uno{x<0}$ should have the same one-dimensional
distribution or, equivalently, that
\begin{equation}
  \Bt(x)-x^2\uno{x<0}\quad{\buildrel {\rm (d)} \over =}\quad\sup_{y\leq x} \big\{\aip(y)-y^2\big\}
  \label{eq:conj}
\end{equation}
for each fixed $x\in\rr$. This distributional identity has actually been proved
rigourously, and its proof is based on the methods we will survey in Sections
\ref{sec:bvalcst} and \ref{sec:appl} (see Theorem \ref{thm:1to2}).

The same heuristic argument works for the other two crossover cases. If we let
$z(0,x)=e^{B(x)}\uno{x\geq0}$, where $B(x)$ is a standard Brownian motion, then
\eqref{eq:1to2approx} and \eqref{eq:zed2} are replaced respectively by
\begin{gather}
  -\log z(t,x) \approx \tfrac1{2t}x^2\uno{x<0}+\tfrac{1}{24}t+\log(\sqrt{2\pi t}) -
  2^{-1/3}t^{1/3}\mathcal{A}_{2\to {\rm BM}}(2^{-1/3}t^{-2/3} x)\\
  \shortintertext{and} z( t,x)=\int_0^\infty dy\,z(0,y; t,x)= \int_0^\infty
  dy\,\tfrac{1}{\sqrt{2\pi t}} e^{ -\tfrac{(x-y)^2}{2t}
    -\tfrac{t}{24}+B(y)+2^{-1/3}t^{1/3}A_t(2^{-1/3}t^{-2/3} (x-y)) },
\end{gather}
and now the same scaling argument allows to conjecture that
\[\mathcal{A}_{2\to {\rm BM}}(x)-x^2\uno{x<0}\quad{\buildrel {\rm (d)} \over
  =}\quad\sup_{y\le x} \big( \aip (y)+\widetilde{B}(x-y)-y^2\big)\] for each fixed
$x\in\rr$, where now $\widetilde B(y)$ is a Brownian motion with diffusion coefficient
2. An analogous argument with $z(0,x)=\uno{x\leq 0}+e^{B(x)}\uno{x\geq0}$ translates into
conjecturing that
\[\mathcal{A}_{1\to {\rm BM}}(x)\quad{\buildrel {\rm (d)} \over =}\quad\sup_{y\in\rr}
\big( \aip (y)+\widetilde{B}(x-y)\uno{y\leq x}-y^2\big)\] for each fixed $x\in\rr$.  As we explained,
these equalities in distribution will only hold in the sense of one-dimensional distributions, i.e. for each fixed
$x$.

The strategy used in the proof of \eqref{eq:conj} is considerably more difficult to
implement for the other two crossover cases (see Section \ref{sec:airy2to1} and the discussion at the end of Section 1.2
in \cite{qr-airy1to2}), and in fact these identities remain conjectures for now. In work
in progress \cite{corwinLiuWang} obtain an improved version of the slow decorrelation
result proved in \cite{corwinFerrariPeche1}, which should allow to prove a general version
of formulas for last passage times in last passage percolation in terms of variational
problems for the Airy$_2$ process. In particular, such a result would give a proof of
these conjectural formulas.

\subsection{Determinantal formulas and extended kernels}
\label{sec:extendedDet}

As we already mentioned, the results we will survey in Sections \ref{sec:bvalcst} and
\ref{sec:appl} are based on alternative Fredholm determinant formulas for the
finite-dimensional distributions of the Airy processes. We will introduce these formulas
in Section \ref{sec:bvalcst}, but before doing that let us explain why the original
extended kernel formulas are natural. We will first do this in a simpler setting, namely
random point processes on finite sets. Then we will explain how similar arguments can be
used to derive the formula \eqref{eq:detform} for the finite-dimensional distributions of the Airy$_2$
process.

\subsubsection{Extended kernels and the Eynard-Mehta Theorem}\label{sec:e-m}

Let $\cx$ be a finite set. A \emph{random point process} on $\cx$ is 
a probability measure on the family $2^{\cx}$ of subsets of $\cx$, which we think of as
point configurations. 
A random point process is called \emph{determinantal} if there exists a $|\cx|\times|\cx|$
matrix $K$ with rows and columns indexed by the elements of $\cx$ such that
\[\rho(A):=\pp\big(\{X\in2^\cx\!:\,A\subseteq X\}\big)=\det(K|_A),\]
where $\pp$ is the probability measure underlying the point process and $K|_A$ is the
submatrix of $K$ indexed by $A$,
\[K_A=\big[K(x,x')\big]_{x,x'\in A}.\]
The function $\rho$ is called the \emph{correlation function} of the process, and $K$ is
called its \emph{correlation kernel}. For more details see \cite{borOxford} and references
therein. The term determinantal was introduced in \cite{bor-ol}.

We are interested in a particular type of random point processes. Let $\cx^1,\dotsc,\cx^n$
be $n$ disjoint finite sets. We consider a point process supported on $kn$-point
configurations with the property that there are exactly $k$ points in each $\cx^i$. The
probability of such a configuration is given as follows: given collections of points
$\{x^j_i\}_{i=1,\dotsc,k}\subseteq\cx^j$ for $j=1,\dotsc,n$, we set
\begin{multline}
  \pp\!\left(\left\{\{x^1_i\}_{i=1,\dotsc,k}\cup\dotsm\cup\{x^n_i\}_{i=1,\dotsc,k}\right\}\right)\\
=Z^{-1}\det\big[\phi_i(x_j^1)\big]_{i,j=1}^k\det\big[W_1(x_i^1,x_j^2)\big]_{i,j=1}^k
\dotsm\det\big[W_{n-1}(x_i^{n-1},x_j^n)\big]_{i,j=1}^k\det\big[\psi_i(x_j^n)\big]_{i,j=1}^k,\label{eq:detstr}
\end{multline}
where the $\phi_i$'s are some functions on $\cx^1$, the $\psi_i$'s are some functions on
$\cx^n$ and the $W_i$'s are matrices with rows indexed by $\cx^i$ and columns indexed
by $\cx^{i+1}$. The normalization constant $Z$ is chosen so that the total mass of the
measure is 1. We are assuming implicitly that the right hand side above is non-negative
for any admissible point configuration.

Write $\Phi$ for the $k\times|\cx^1|$ matrix with $k$ rows and columns indexed by elements of $\cx^1$
which is defined by $\Phi_{i,x}=\phi_i(x)$ for $1\leq i\leq k$ and $x\in\cx^1$. Similarly,
write $\Psi$ for the $|\cx^n|\times k$ matrix with $k$ columns and rows indexed by $\cx^k$
which is defined by $\Psi_{x,i}=\psi_i(x)$ for $1\leq i\leq k$ and
$x\in\cx^n$. Furthermore, define the $k\times k$ matrix
\[M=\Phi W_1\dotsm W_{n-1}\Psi.\] We will assume that $\det(M)\neq0$. Under this
assumption it can be shown (see e.g. \cite{borodinRains}) that the normalization constant
$Z$ in \eqref{eq:detstr} equals $\det(M)$.

The Eynard-Mehta Theorem states that a random point processes defined as in
\eqref{eq:detstr} is determinantal. Moreover, the theorem gives an explicit formula for
the correlation kernel. The precise statement is the following:

\begin{thm}[\cite{eynardMehta}]\label{thm:eynardMehta}
  The random point processes defined by \eqref{eq:detstr} is determinantal. Its
  correlation kernel is the block matrix $K$ with $n\times n$ blocks, such that the
  $(i,j)$-block has rows indexed by $\cx^i$ and columns indexed by $\cx^j$, and is given by
  \[K_{i,j}=W_i\dotsm W_{n-1}\Psi M^{-1}\Phi W_1\dotsm W_{j-1}-W_i\dotsm W_{j-1}.\]
\end{thm}

For a simple proof of this result see \cite{borodinRains}. Remarkably, the inverse
$M^{-1}$ can be computed, or at least approximated, in many cases of interest.

The connection with the models we have discussed so far is through certain families of
non-intersecting paths. The Airy$_2$ process can be obtained directly as a limit of the
top line of several different families of non-intesecting paths, one of which is presented
in Section \ref{sec:dysonBM} (for some others see \cite{johanssonRMandDetPr}). For the
other Airy processes presented in Section \ref{sec:spatialfluc} the connection with
non-intersecting paths is less immediate (see for instance the discussion preceding Lemma
3.4 in \cite{borFerPrahSasam}), but in any case enough of the above structure remains, and
the proofs still rely crucially on a version of Theorem \ref{thm:eynardMehta}. On the
other hand, we may think of the $kn$-point configurations where the measure defined in
\eqref{eq:detstr} is supported as defining a family of $k$ (in principle not necessarily
non-intersecting) paths. For example, the first path would be expressed by
$(x^1_1,x^2_1,\dotsc,x^n_1)$. It turns out, as we will see below, that probability
measures on families of $k$ non-intersecting paths on $\cx^1\cup\dotsm\cup\cx^n$ are
naturally given by expressions like \eqref{eq:detstr}, and hence have a determinantal
structure. If the sets $\cx^i$ are endowed with some total order and we assume that our
non-intersecting paths are arranged so that $(x^1_1,x^2_1,\dotsc,x^n_1)$ is the top path,
then one can prove (see e.g. \cite{johansson}) that
\begin{equation}
\pp\!\left(x^1_1\leq z_1,\dotsc,x^n_1\leq z_n\right)=\det(I-PKP),\label{eq:fd}
\end{equation}
where $z_i\in\cx^i$, $K$ is the correlation kernel given by the Eynard-Mehta Theorem and
$P$ is block-diagonal matrix with $n$ diagonal blocks defined so that, for $i=1,\dotsc,n$,
$P_{i,i}$ has rows and columns indexed by $\cx^i$ and is given by
$(P_{i,i}v)_j=\uno{x^i_j>z_i}v_j$. This should be compared with an expression like
\eqref{eq:detform}.

If we go back to thinking about these paths as defining a random point process, then they
are given by a measure on $kn$-point configurations on
$\cx^1\cup\dotsm\cup\cx^n$. Therefore, if the process is determinantal, its correlation
kernel necessarily has to be a matrix with rows and colums indexed by
$\cx^1\cup\dotsm\cup\cx^n$. The Eynard-Mehta Theorem implies moreover that the correlation
matrix is partitioned naturally into $n\times n$ blocks, with the $(i,j)$ block having rows
indexed by $\cx^i$ and colums indexed by $\cx^j$. To see how this structure relates with
the extended kernels introduced in Section \ref{sec:spatialfluc} for the Airy processes,
we make the following observation. An operator $T$ acting on
$L^2(\{t_1,\dotsc,t_n\}\times\rr)$ can be regarded as an operator-valued matrix
$\big(T_{i,j}\big)_{i,j=1,\dots,n}$ with entries $T_{i,j}$ (acting on $L^2(\rr)$), which
acts on $f\in L^2(\rr)^{n}$ as $(Tf)_i=\sum_{j=1}^n T_{i,j}f_j$ (more precisely, we are
using the fact that $L^2(\{t_1,\dotsc,t_n\}\times\rr)$ and
$\bigoplus_{t\in\{t_1,\dotsc,t_n\}}L^2(\rr)$ are isomorphic as Hilbert spaces). Hence an
extended kernel formula like \eqref{eq:detform} can be thought of as the determinant of an
$n\times n$ matrix whose entries are operators on $L^2(\rr)$. Similarly, we may think of
\eqref{eq:fd} as the determinant of an $n\times n$ matrix whose $(i,j)$ entry maps
$L^2(\cx^i)$ to $L^2(\cx^j)$. Since the Airy processes live on the real line instead of
finite sets, this latter spaces are replaced by $L^2(\rr)$.

In the next section we will explain how these ideas can be used to derive the extended
kernel formula for Airy$_2$.

\subsubsection{Derivation of the \texorpdfstring{Airy$_2$}{Airy2} process from Dyson
  Brownian motion}\label{sec:dysonBM}

The original derivation of the Airy$_2$ process was done in \cite{prahoferSpohn} using
quantum statistical mechanical arguments, while Johansson's proof of
Theorem \ref{thm:Airy2LPP} relies crucially on the connection between LPP and the
Robinson-Schensted-Knuth algorithm, which provides a family of discrete non-intersecting
paths, the top line of which converges to $\aip$. We will briefly explain the derivation
of the Fredholm determinant formula for the Airy$_2$ process using a different model, the Dyson
Brownian motion. We refer the reader to
\cite{tracyWidomNIBrExc,tracyWidomDysonBM,johansson,johanssonRMandDetPr} for more details
on the derivation of Airy$_2$ from non-intesecting paths.

Consider the evolving eigenvalues of an $N\!\times\!N$ GUE matrix with each (algebraically
independent) entry diffusing according to a stationary Ornstein-Uhlenbeck process. We
write the eigenvalues at time $t$ as $\lambda^{N}(t) =
(\lambda^N_1(t),\dots,\lambda^N_N(t))$ so that $\lambda_i(t)$ decreases in $i$. This
eigenvalue diffusion, called the \emph{stationary GUE Dyson Brownian motion}, can be
written as the solution of a certain $N$-dimensional SDE, and it can be shown that it is
stationary, with distribution given by the eigenvalue distribution of an $N\times N$ GUE
matrix. Moreover, the paths followed by the $N$ eigenvalues almost surely form an ensemble
of non-intersecting curves.

Suppose we look at this eigenvalue diffusion at times $t_1<\dotsm<t_n$, and we condition
the $N$ paths to be pairwise non-intersecting. To investigate the transitions
between $t_m$ and $t_{m+1}$, suppose we condition
this eigenvalue diffusion to start at time $t_m$ at $\lambda^N_i(t_m)=x_i$ for some fixed
$x_1<\dotsm<x_N$, and we also fix destination points $y_1<\dotsm<y_N$. Let $p_t(x,y)$ be the
transition probability density of a (one-dimensional) Ornstein-Uhlenbeck process from $x$
at time 0 to $y$ at time $t$. Then, in this setting, the Karlin-McGregor Theorem
\cite{karlinMcGregor} implies that the transition probability density for these $N$
non-intersecting paths to end at the prescribed destination points $y_1,\dotsc,y_N$ is
given by a constant times
\begin{equation}\label{eq:kmcg}
  \det\!\big[p_{t_{m+1}-t_m}(x_i,y_j)\big]_{i,j=1}^N.
\end{equation}
The transition function $p_{t_{m+1}-t_m}$ corresponds then to the matrix $W_m$ in
\eqref{eq:detstr}. Of course, since our paths take values in $\rr$ now, we no longer have
a matrix, but the Eynard-Mehta Theorem still holds in this setting (see
e.g. \cite{tracyWidomDysonBM}). The functions $\phi_i$ and $\psi_i$ in \eqref{eq:detstr}
are related to the (stationary) marginals for $\lambda^N_t$, and in this case are equal
(due to stationarity) and expressed simply in terms of Hermite polynomials. The result
after further computations and using \eqref{eq:fd} is the following
\cite{tracyWidomDysonBM}: given $x_1,\dotsc,x_n\in\rr$,
\begin{equation}\label{GUEdpp}
  \pp\!\left(\lambda^N_1(t_1)\leq x_1,\dotsc,\lambda^N_1(t_n)\leq
    x_n\right)=\det\!\big(I-{\rm f}K^{\rm ext}_{{\rm Hrm},N}{\rm f}\big)_{L^2(\{t_1,\dots,t_n\}\times\rr)},
\end{equation}
where ${\rm f}$ is defined as in \eqref{eq:deff} and $K^{\rm ext}_{{\rm Hrm},N}$ is the
\emph{extended Hermite kernel}
\[K^{\rm ext}_{{\rm Hrm},N}(s,x;t,y)=
\begin{dcases*}
  \sum_{k=0}^{N-1}e^{k(s-t)}\varphi_k(x)\varphi_k(y) & if $s\geq t$,\\
  -\sum_{k=N}^\infty e^{k(s-t)}\varphi_k(x)\varphi_k(y) & if $s<t$,
\end{dcases*}\] and where $\varphi_k(x)=e^{-x^2/2}p_k(x)$ with $p_k$ the $k$-th normalized
Hermite polynomial (so that $\|\varphi_k\|_2=1$). Note that $\lambda^N_1$ is the top line
of our family of non-intersecting paths, so this probability is the same as the
probability that all paths stay below the $x_i$'s. Note also the similarity between this
formula and the formula \eqref{eq:detform} for $\aip$. We remark that the scaling
of the eigenvalues appearing in the last formula differs by a factor of $\sqrt{N}$ with
the one introduced in Section \ref{sec:tw}; the present choice is the one that is
naturally associated with the operator $D$ introduced next.

The kernel $K^{\rm ext}_{{\rm Hrm},N}$ has a nice algebraic
structure. Writing \[D=-\tfrac12\big(\Delta-x^2+1\big),\]
(i.e. $Df(x)=-\frac12(f''(x)-(x^2-1)f(x)))$, the \emph{harmonic oscillator functions}
$\varphi_k$ satisfy $D\varphi_k=k\varphi_k$, and moreover $\big\{\varphi_k\}_{k\geq0}$
forms a complete orthonormal basis of $L^2(\rr)$. Define the \emph{Hermite kernel} as
\[K_{{\rm Hrm},N}(x,y)=\sum_{k=0}^{N-1}\varphi_k(x)\varphi_k(y),\] which is then just the
projection onto span$\{\varphi_0,\dotsc,\varphi_{N-1}\}$. Then the following formula
holds:
\[K^{\rm ext}_{{\rm Hrm},N}(s,x;t,y)=
\begin{dcases*}
  e^{(s-t)D}K_{{\rm Hrm},N}(x,y) & if $s\geq t$,\\
  -e^{(s-t)D}(I-K_{{\rm Hrm},N})(x,y) & if $s<t$.
\end{dcases*}\]

Now introduce the rescaled process
\[\widetilde\lambda^N_i(t)=\sqrt{2}N^{1/6}\big(\lambda^N_i(N^{-1/3}t)-\sqrt{2N}\big).\]
Changing variables $x\mapsto \frac1{\sqrt{2}N^{1/6}}x+\sqrt{2N}$, $y\mapsto
\frac1{\sqrt{2}N^{1/6}}y+\sqrt{2N}$ in the kernel accordingly, a calculation gives
\[\pp\!\left(\widetilde\lambda^N_1(t_1)\leq x_1,\dotsc,\widetilde\lambda^N_1(t_n)\leq
  x_n\right)=\det\!\big(I-{\rm f}\widetilde K^{\rm ext}_{{\rm Hrm},N}{\rm
  f}\big)_{L^2(\{t_1,\dots,t_n\}\times\rr)}\] with
\[\widetilde K^{\rm ext}_{{\rm Hrm},N}(s,x;t,y)=
\begin{dcases*}
  e^{(s-t)H_N}\widetilde K_{{\rm Hrm},N}(x,y) & if $s\geq t$,\\
  -e^{(s-t)H_N}(I-\widetilde K_{{\rm Hrm},N})(x,y) & if $s<t$,
\end{dcases*}\] where $\widetilde K_{{\rm Hrm},N}(x,y)=\frac1{\sqrt{2}N^{1/6}}K_{{\rm
    Hrm},N}\!\left(\frac{x}{\sqrt{2}N^{1/6}}+\sqrt{2N},
  \frac{y}{\sqrt{2}N^{1/6}}+\sqrt{2N}\right)$ and the operator
$H_N=-\Delta+x+\frac{x^2}{2N^{2/3}}$.

The above rescaling corresponds to focusing in on the top curves of the Dyson Brownian
motion. It is known that in the limit $N\to\infty$, $\widetilde K_{{\rm Hrm},N}$ converges
to the Airy kernel $\K$, while it is clear that $H_N$ converges to the \emph{Airy
  Hamiltonian} $H$:
\begin{equation}
  \label{eq:airyHamilt}
  H=-\Delta+x
\end{equation}
(i.e. $Hf(x)=-f''(x)+xf(x)$). Putting aside precise convergence issues, the result is that
\begin{equation}
\lim_{N\to\infty}\widetilde K^{\rm ext}_{{\rm Hrm},N}(s,x;t,y)=
\begin{dcases*}
  e^{(s-t)H}\K(x,y) & if $s\geq t$,\\
  -e^{(s-t)H}(I-\K)(x,y) & if $s<t$.
\end{dcases*}\label{eq:limHrm}
\end{equation}
The obvious question at this point is what is the relationship between this limit and the
extended Airy kernel \eqref{eq:extAiry}. It turns out that they are the same. This is a
consequence of the following remark, which implies that the nice structure we saw in
$\widetilde K_{{\rm Hrm},N}$ survives in the limit:

\begin{rem}\label{airyrem}
  The shifted Airy functions $\phi_\lambda(x)=\Ai(x-\lambda)$ are the generalized
  eigenfunctions of the Airy Hamiltonian, as $H\phi_\lambda=\lambda\phi_\lambda$ (we say
  generalized because $\phi_\lambda\notin L^2(\rr)$). The Airy kernel $K_{\Ai}$ is the
  projection of $H$ onto its negative generalized eigenspace. This is seen by observing
  that if we define the operator $A$ to be the \emph{Airy transform},
  $Af(x):=\int_{-\infty}^\infty dz\Ai(x-z)f(z)$, then $K_{\Ai}=A\bar P_0A^*$, where $\bar
  P_0f(x)=\uno{x<0}f(x)$.
\end{rem}

In particular, $e^{tH}$ is defined spectrally. Formally, its integral kernel is given by
$e^{tH}(x,y)=\int_{-\infty}^\infty
d\lambda\,e^{-t\lambda}\Ai(x+\lambda)\Ai(y+\lambda)$. The integral converges when $t<0$ by
the decay properties of the Airy function, but it diverges when $t>0$ (it can be
interpreted as $\delta_{x=y}$ when $t=0$). Nevertheless, in our formulas $e^{tH}$ will
always appear after $\K$ when $t>0$. This has the effect of restricting the integral to
$\lambda>0$, which converges because the Airy function is bounded.

As a consequence of the above discussion we obtain the following result, see \cite{tracyWidomDysonBM}:

\begin{thm}
  \[\widetilde\lambda^N_1(t)\xrightarrow[N\to\infty]{}\aip(t)\]
  in the sense of convergence of finite-dimensional ditributions.
\end{thm}

The extended kernels which define the other Airy processes do not have exactly the same
structure. One reason behind this is that, apart from $\aip$ and $\aipo$, the other
processes are not stationary. 
The other reason is that, as we mentioned, in some cases (for example Airy$_1$) the
processes are obtained after further limiting procedures, which destroy part of the
structure. Nevertheless, as we will see later, enough of this structure remains for our
purposes.

We have wandered a bit far from the main subject of this survey in the hope that the
reader will get a feeling about why extended kernels appear naturally for our
processes. In the rest of this article we will deal with formulas which are given as
Fredholm determinants of certain operators acting on $L^2(\rr)$, as opposed
$L^2(\{t_1,\dots,t_n\}\times\rr)$. In light of the above discussion, it is slightly
surprising that such formulas should exist.

\section{Fredholm Determinants}
\label{sec:fred}

If $K$ is an integral operator acting on $H= L^2(X,d\mu)$ through its kernel
\begin{equation}\label{intop}(Kf)(x) = \int_{X} K(x,y)f(y) d\mu(y),\end{equation} we
define the {\it Fredholm determinant} by
\begin{equation}\label{FD}
  \det(I+\lambda K) = 1+\sum_{n=1}^{\infty} \frac{\lambda^n}{n!} \int_{X} \cdots \int_{X} \det\left[K(x_i,x_j)\right]_{i,j=1}^{n} d\mu(x_1)\cdots d\mu(x_n).
\end{equation}
If $|K(x,y)|\le B$ for all $x,y$, and $\mu$ is a finite measure, the {\it Fredholm series}
\eqref{FD} converges by Hadamard's inequality,
\begin{equation*}
  \left|\det( C_1,\ldots, C_n)\right|\le  \|C_1\|\cdots \|C_n\|
\end{equation*}
where $\|C_i\|$ denotes the Euclidean length of the column vector $C_i$, since the length
of the column vector in $\left[K(x_i,x_j)\right]_{i,j=1}^{n}$ is bounded by $Bn^{1/2}$,
and hence the $n$-th summand in (\ref{FD}) is bounded by $\tfrac{\lambda^n}{n!} B^n
n^{n/2}$.

If one is not familiar with the definition \eqref{FD} one might even wonder what it has to
do with determinants.  Take a matrix $K=\left[K_{ij}\right]_{i,j=1}^{d}$, $d<\infty$, and
consider the $d\times d$ determinant $\det(I+\lambda K)$.  Clearly it is a polynomial of
degree $d$ in $\lambda$, $\sum_{n=0}^d a_n\lambda^n$, and its coefficients are given by
the rule $a_n =\tfrac{1}{n!}\partial_\lambda^n \det(I+\lambda K)\!\mid_{\lambda=0} $.  To
compute this, use the rule for differentiating determinants,
\begin{equation*}
  \partial_\lambda \det( C_1,\ldots, C_d)=\sum_{n=1}^d \det( C_1,\ldots, \partial_\lambda C_n, \ldots, C_d)
\end{equation*}
and the fact that, in our particular case, $C_n(\lambda)= e_n+ \lambda K_{\cdot,n}$ is
linear in $\lambda$ and $C_n(0)= e_n$, the $n$-th unit vector.  The result is
\begin{multline}
  \det(I+\lambda K)=1+ \lambda\sum_{1\le i\le d} K_{ii} + \lambda^2\sum_{1\le i<j\le d} \det \left[ \begin{array}{cc}
      K_{ii} & K_{ij} \\
      K_{ji} & K_{jj}  \end{array} \right] \\
  +\lambda^3\sum_{1\le i<j<k\le d}
  \det  \left[\begin{array}{ccc}
      K_{ii} & K_{ij} & K_{ik} \\
      K_{ji} & K_{jj} & K_{jk} \\
      K_{ki} & K_{kj} & K_{kk} \end{array} \right]+\cdots+ \lambda^d\det K.
\end{multline}
Replacing the ordered sums with unordered sums gives a factor ${1}/{n!}$, and setting
$\lambda=1$ we can see that this is a special case of \eqref{FD}. Von Koch's idea
\cite{vonKoch} was that this formula for the determinant was the natural one to extend to
$d=\infty$.  Fredholm replaced the integral operator \eqref{intop} on $L^2([0,1], dx)$ by
its discretization $[ \tfrac1{n}K(\tfrac{i}{n}, \tfrac{j}{n}) ]_{i,j=1}^n$ to obtain
\eqref{FD}, which he then used to characterize the solvability of the integral equation
$(I+K) u=f$ via the non-vanishing of the determinant of $I+K$.

One can of course imagine other, more intuitive definitions of the determinant.  Perhaps
\begin{equation}\label{proddet}
  \det(I+K) = \prod_n (1+\lambda_n)
\end{equation}
where $\lambda_n$ are the eigenvalues of $K$, counted with multiplicity.  Or
\begin{equation}\label{DSdet}
  \det(I+\lambda K) = e^{\tr \log(1+\lambda A)}
\end{equation}
with the trace
\begin{equation}
\tr K = \int d\mu(x)\,K(x,x).\label{eq:trace}
\end{equation}
Of course, these definitions require some smallness condition on $K$, but at least they
make apparent the important fact that the determinant is invariant under conjugation
$\det(I+ M^{-1}KM) =\det(I+K)$, or
\begin{equation}\label{eq:cyclic}
  \det(I+ K_1K_2) = \det(I+K_2K_1),
\end{equation}
(usually referred to as the cyclic property of determinants) as well as the formula
\begin{equation}\label{derivativeofadet}
  \partial_\beta \det(I+K(\beta)) = \det(I+K(\beta))\tr( (I+K(\beta))^{-1}\partial_\beta K(\beta))
\end{equation}
for $K(\beta)$ depending smoothly on a parameter $\beta$.

A more modern way to write \eqref{FD} is
\begin{equation}\label{Sdet}
  \det(I+\lambda K) =\sum_{n=0}^\infty \lambda^n \tr \Lambda^n(K)
\end{equation}
where $\Lambda^n(K)$ denotes the action of the tensor product $A\otimes \cdots\otimes A$
on the antisymmetric subspace of $H\otimes \cdots\otimes H$.  If $P_n$ denotes the
projection onto that subspace and $C_n = P_n \Lambda^n(K) P_n$ then
\begin{eqnarray*}
  &&C_n(f_1\otimes\cdots\otimes f_n) =  \frac1{n!}\sum_{\sigma \in S_n} {\mathop{\rm sgn}}(\sigma) Af_{\sigma(1)} \otimes\cdots\otimes Af_{\sigma(n)} \\
  &&= \frac1{n!}\sum_{\sigma \in S_n} {\mathop{\rm sgn}}(\sigma)\int\!\cdots\!\int
  d\mu(y_1)\cdots d\mu(y_n)\,K(x_1,y_{\sigma(1)})\cdots K(x_n,y_{\sigma(n)})
  f_1(y_1)\cdots f_n(y_n) 
\end{eqnarray*}
which shows that $C_n$ is an integral operator with kernel
$\det\left[K(x_i,x_j)\right]_{i,j=1}^{n}$ and hence \eqref{Sdet} is just a slick way to
write \eqref{FD}. The advantage of \eqref{FD} is that it can be used directly to define
the Fredholm determinant for operators on a general separable Hilbert space, but we will
not need this point of view here (see \cite{simon} for more details).

The natural notion of smallness for Fredholm determinants turns out to be the trace norm
on operators
\[\|K\|_1 := \tr |K|,\]
where $|K|=\sqrt{K^*K}$ is the unique positive square root of the operator $K^*K$.  A
(necessarily compact) operator with finite trace norm is called {\it trace class}.  Using
the Parseval relation, one can check that for such operators the trace can be defined as
\begin{equation}
  \tr K=\sum_{n=1}^{\infty}\langle e_n,K e_n\rangle,
\end{equation}
as it is basis independent. This works for operators on any separable Hilbert space, and
in our setting it can be shown that this definition of trace coincides with \eqref{eq:trace}
for $K$ of trace class. The \emph{Hilbert-Schmidt norm} $\|K\|_2 =\sqrt{\tr(|K|^2)}$ is
easier to compute,
\[\|K\|_2=\left(\int dx\,dy\,|K(x,y)|^2\right)^{1/2},\] and the relation between these
norms and the more common operator norm $\|K\|_\mathrm{op}$ is
\begin{equation}\|K\|_\mathrm{op}\leq\|K\|_2\leq\|K\|_1,\end{equation} as well as
\[\|K_1K_2\|_1\leq\|K_1\|_2\|K_2\|_2,\quad
\|AK\|_1\leq\|A\|_\mathrm{op}\|K\|_1,\qand\|AK\|_2\leq\|A\|_\mathrm{op}\,\|K\|_2,\] all of
which can be checked easily.  Of course, in the latter two $A$ need not be compact.  The
reason the trace norm is so useful is

\begin{lem}\label{lem:fredholm}\mbox{}
  \begin{enumerate}[label=\arabic*.]
  \item (Lidskii's Theorem) If $K$ is trace class then $ \tr K = \sum_n\lambda_n$, where
    $\lambda_n$ are the eigenvalues of $K$. It follows that the three definitions
    \eqref{FD}, \eqref{proddet} and \eqref{DSdet} are equivalent.
  \item $A\mapsto\det(I+A)$ is continuous in trace norm. Explicitly,
    \begin{equation}\label{contintr}\left|\det(I+K_1)-\det(I+K_2)\right|\leq\|K_1-K_2\|_{1}\exp(\|K_1\|_1+\|K_2\|_1+1).\end{equation}
  \end{enumerate}
\end{lem}

Lidskii's theorem is non-trivial and its proofs use heavy function theory, but
\eqref{contintr} can be explained easily.  Let $f(z) = \det(I+ \tfrac12 (K_1+K_2) +
z(K_1-K_2))$, so that the left hand side of \eqref{contintr} is $|f(\tfrac12)
-f(-\tfrac12)|\le \sup_{-1/2\le t\le 1/2} |f'(t)|$.  Cauchy's integral formula $f'(z)
=\tfrac1{2\pi i} \oint \frac{f(z')}{z'-z} dz'$ shows that $ \sup_{-1/2\le t\le 1/2}
|f'(t)|\le \tfrac1{R} \sup_{|z|\le R+\tfrac12} |f(z)| $.  The eigenvalues of
$\Lambda^n(K)$ are $\lambda_{i_1}\cdots \lambda_{i_n}$, $i_1<\cdots<i_n$, so $ \tr
\Lambda^n(K)= \sum_{i_1<\cdots<i_n} \lambda_{i_1}\cdots \lambda_{i_n} $ and hence
$\left|\tr \Lambda^n(K)\right| \le \tfrac1{n!} \|K\|_1^n $, which implies
\begin{equation}|\det(I+\lambda K)| \le e^{\lambda \|K\|_1}.
\end{equation}
Therefore $\sup_{|z|\le R+\tfrac12} |f(z)| \le \exp(\tfrac12 \|K_1+K_2\|_1+(R+\tfrac12)
\|K_1-K_2\|_1)$ and taking $R=\| K_1-K_2\|_1^{-1}$ gives \eqref{contintr}.

\vspace{2pt}
\begin{exs}\mbox{}\\
  {\it 1. (Gaussian distribution)} A trivial example is $K(x,y) =
  \frac{e^{-x^2/2t}}{\sqrt{2\pi t}}$.  The operator is rank one, so if $P_s$ is the
  orthogonal projection from ${L^2(\rr)} \to L^2(s,\infty)$ then by \eqref{DSdet} we have
  \begin{equation}
    \det( I - P_s K P_s)_{L^2(\rr,dx)} = 1- \tr P_s K P_s = \int_{-\infty}^sdx\, \frac{e^{-x^2/2t}}{\sqrt{2\pi t}}.
  \end{equation}
  Of course, the Gaussian here could be replaced by an arbitrary density.\\[6pt]
  \noindent {\it 2. (GUE)} Consider the Airy kernel $K_{\rm Ai} (x,y)=\int_0^\infty
  dt\Ai(x+t)\Ai(y+t)$ and let $\Ai_t(x)=\Ai(x+t)$ and $H=-\partial_x^2 +x$.  Then $H{\rm
    Ai}_t = -t {\rm Ai}_t$, the ${\rm Ai}_t$, $t\in\rr$ are generalized eigenfunctions of
  $H$, and $K_{\rm Ai}$ is the orthogonal projection onto the negative eigenspace of $H$
  (see Remark \ref{airyrem}). Using ${\rm Ai}''(x)=x{\rm Ai}(x)$, we have $\partial_t
  \frac{ {\rm Ai}(x+t) {\rm Ai}'(y+t) - {\rm Ai}'(x+t) {\rm Ai}(y+t)}{y-x} ={\rm Ai}(x+t)
  {\rm Ai}(y+t) $, which yields the Christoffel-Darboux formula 
  \begin{equation}
    K_{\rm Ai} (x,y) = \frac{ {\rm Ai}'(x+t) {\rm Ai}(y+t) - {\rm Ai}(x+t) {\rm Ai}'(y+t)}{y-x}.
  \end{equation}
  To show $P_sK_{\rm Ai}P_s$ is trace class, write $K_{\rm Ai} =B_0P_0B_0$ where
  \begin{equation}
    \label{bee}
    B_0(x,y) ={\rm Ai} (x+y).
  \end{equation}
  Then use $\|K_1K_2\|_1\le \|K_1\|_2\|K_2\|_2$ to get
  \begin{equation}\label{atest}
    \|P_sK_{\rm Ai}P_s\|_1 \le\|P_sB_0P_0\|_2^2\le \int_0^\infty\int_s^\infty {\rm Ai}^2(x+y) dxdy,
  \end{equation}
  which is finite by the following well-known estimates for the Airy function (see
  (10.4.59-60) in \cite{abrSteg}):
  \begin{equation}
    |\!\Ai(x)|\leq Ce^{-\frac23x^{3/2}}\quad\text{for $x>0$},\qquad
    |\!\Ai(x)|\leq C\quad\text{ for $x\leq0$}.\label{eq:airybd}
  \end{equation}
  The GUE Tracy-Widom distribution is given by
  \begin{equation}\label{TWGUE}
    F_{\rm GUE}(s) = \det(I-P_sK_{\rm Ai}P_s)_{L^2(\rr,dx)}.
  \end{equation}
  On the face of it, it is not so obvious why such an expression would define a
  probability distribution function.  From \eqref{atest} it is clear that $\lim_{s\to
    \infty} \det(I-P_sK_{\rm Ai}P_s)=1$.  Since $P_s\K P_s$ is a composition of
  projections, its eigenvalues satisfy $1\ge \lambda_1(s)\ge \lambda_2(s) \ge\cdots\ge 0$.  Recall
  the min-max characterization of eigenvalues
  \begin{equation}
    \lambda_k(s) = \max_{\dim U=k} \min_{f\in U} \frac{ \langle f, P_sK_{\rm Ai}P_s f\rangle}{\langle f,f\rangle},
  \end{equation}
  from which it is apparent that $\lambda_i(s)$ is non-decreasing as $s$ decreases, and
  hence $\det(I-P_sK_{\rm Ai}P_s)=\prod_i(1-\lambda_i(s))$ is non-increasing with
  decreasing $s$.  In fact, $\lambda_1(s)\nearrow 1$ as $s\searrow -\infty$ since if $f$
  is in the negative eigenspace of $H$, $\langle P_s f,K_{\rm Ai}P_s f\rangle \to \langle
  f, K_{\rm Ai} f\rangle= \langle f, f\rangle$. This shows that
  $\det(I-P_sK_{\rm Ai}P_s)\searrow 0$ as $s\searrow -\infty$ (for an asymptotic expansion
  of $F_{\rm GUE}(s)$ as $s\searrow-\infty$ see \cite{baikBuckDiF}).\\[6pt]
  \noindent {\it 3. (GOE)} $F_{\rm GOE}(s) = \det(I-P_sB_0P_s)_{L^2(\rr,dx)}$ where
  $B_0(x,y)$ is as in \eqref{bee}. The key to show that $B_0$ is trace class in this case
  is the identity
  \begin{equation}
    \int_{-\infty}^\infty dx\Ai(a+x)\Ai(b-x)=2^{-1/3}\Ai(2^{-1/3}(a+b))\label{eq:airyIdentity}
  \end{equation}
  (see, for example, (3.108) in \cite{valleeSoares}). One defines
  $G_1(x,z)=2^{1/3}\!\Ai(2^{1/3}x+z)e^{z}$ and $G_2(z,y)=e^{-z}\Ai(2^{1/3}y-z)$ and notes
  that $P_sB_0P_s=(P_sG_1)(G_2P_s)$. Then \eqref{eq:airybd} allows to show that each of
  the last two factors has finite Hilbert-Schmidt norm, yielding that $P_sB_0P_s$ is
  trace class.\\[6 pt]

  {\it 4. (Airy$_1$ process)} Recall the Fredholm determinant formula \eqref{eq:detAiry1}
  for the finite-dimensional distributions of the Airy$_1$ process. It turns out that the
  kernel ${\rm f}K^{\rm ext}_1{\rm f}$ inside the determinant is not trace
  class, basically because the heat kernel is not even Hilbert-Schmidt on
  $L^2([s,\infty))$ for $s\in\rr$. Nevertheless, the series \eqref{FD} defining the
  Fredholm determinant is finite in this case, because one can conjugate the kernel ${\rm
    f}K^{\rm ext}_1{\rm f}$ to something which can be proved to be trace class
  (see \cite{bfp}).
\end{exs}

The situation in the last example, where the natural expression for a kernel defines an
operator which is not trace class, but which is conjugate to a trace class operator,
arises often. Here by conjugacy we mean the following: two operators $K$ and $\wt K$ are
conjugate if there exists some invertible linear mapping $U$ acting on measurable
functions on $X$ such that $K=U\wt KU^{-1}$. Observe that such a pair of operators have
the same Fredholm series expansion \eqref{FD}, i.e. $\det(I+K)=\det(I+\wt K)$. This allows
to extend the manipulations on Fredholm determinants to operators which are
conjugate to trace class operators, provided that one is careful in keeping track of the
needed conjugations.
  
The reason we start with the Fredholm expansion \eqref{FD} is that this is the way the
determinant usually arises from combinatorial expressions. Sometimes the kernels are not
trace class, but this should not bother us so much as long as some version of the formal
expression can be shown to converge, for instance as in Example 4 above.  Often, it is
genuinely difficult to show that the resulting expressions define a probability
distribution, and we only know it because they arose this way.

\section{Boundary Value Kernels and Continuum Statistics of Airy Processes}
\label{sec:bvalcst}

\subsection{Boundary value kernel formulas for finite-dimensional distributions}\label{sec:bdVlK}

Recall the formula \eqref{eq:detform} for the finite-dimensional distributions of the
Airy$_2$ process. It is given in terms of the Fredholm determinant of what we call an
\emph{extended kernel}, that is, (the kernel of) an operator acting on the ``extended
space'' $L^2(\{t_1,\dotsc,t_n\}\times\rr)$. Although such formulas have been very useful
in the study of models in the KPZ class, they suffer from two problems. First, if one
wants to take the number $n$ of times $t_i$ to infinity, a big difficulty appears in the
fact that these formulas involve Fredholm determinants on the Hilbert space
$L^2(\{t_1,\dotsc,t_n\}\times\rr)$, and thus the space itself is changing as $n$
grows. Second, these formulas are useful for computing long range properties of the
processes (for instance an asymptotic expansion of the covariance of $\aip(s)$ and
$\aip(t)$ as $|t-s|\to\infty$, see \cite{widomAiry2}), but are not suitable for studying short
range properties such as regularity of the sample paths.

The second type of Fredholm determinant formula, which is the one we will use for most of
the rest of this article, was actually introduced as the original definition of the
Airy$_2$ process by \citet{prahoferSpohn}. It is given as follows: for $t_1<\dotsm<t_n$
and $x_1,\dotsc,x_n\in\rr$,
\begin{multline}
  \label{eq:airyfd}
  \pp\!\left(\aip(t_1)\leq x_1,\dotsc,\aip(t_n)\leq x_n\right)\\
  =\det\!\left(I-K_{\Ai}+\bar P_{x_1}e^{(t_1-t_2)H}\bar P_{x_2}e^{(t_2-t_3)H}\dotsm \bar
    P_{x_n}e^{(t_n-t_1)H}K_{\Ai}\right)_{L^2(\rr)},
\end{multline}
where $K_{\Ai}$ is Airy kernel \eqref{eq:airyKernel}, $H$ is the Airy Hamiltonian
\eqref{eq:airyHamilt} and $\bar P_a$ denotes the projection onto the interval
$(-\infty,a]$:
\[\bar P_af(x)=\uno{x\leq a}f(x).\]
Note that the Fredholm determinant is now computed on the Hilbert space $L^2(\rr)$ instead
of $L^2(\{t_1,\dotsc,t_n\}\times\rr)$, which makes taking $n\to\infty$ at least
feasible. Note also that the time increments $t_i-t_{i+1}$ appear explicitly in the
formula, which explains why this formula will be more suitable for the study of short
range properties. Another advantage of this formula is that it makes apparent that $\aip$
is a stationary process.

The equivalence of \eqref{eq:detform} and \eqref{eq:airyfd} was derived formally in
\cite{prahoferSpohn} and \cite{prolhacSpohn}. The proof in \cite{prolhacSpohn} is based in
the following idea. As we explained in Section \ref{sec:e-m}, the extended kernel formula
\eqref{eq:detform} can be thought of as the determinant of an $n\times n$ matrix whose
entries are operators acting on $L^2(\rr)$. By rewriting this operator as a sum of an
upper-triangular part and lower-triangular part and using algebraic properties of the
determinant and the algebraic relationships between the different entries of this matrix,
\cite{prahoferSpohn} showed that the determinant equals the determinant of an
operator-valued matrix $I+G$ such that only the first column of $G$ is non-zero. Therefore
$\det(I+G)_{L^2(\{t_1,\dotsc,t_n\}\times\rr)}=\det(I+G_{1,1})_{L^2(\rr)}$ (to see this
simply pretend the operators in the determinants are matrices), and an explicit
calculation of $G_{1,1}$ yields \eqref{eq:airyfd}.

The argument given in \cite{prolhacSpohn} which we just sketched is almost a complete
proof. There are nevertheless some subtleties. For example, it is not a priori obvious that
for $s,t>0$, $e^{-sH}$ can be applied to the image of $\bar P_ae^{-tH}$. Moreover, in
order to manipulate Fredholm determinants one needs to check that certain analytical
conditions are satisfied (see Section \ref{sec:fred}). The technical details are discussed
in \cite{quastelRemAiry1}, which in fact shows that a formula analogous to
\eqref{eq:airyfd} holds for the Airy$_1$ process as well. It is given as follows: for
$t_1<\dotsm<t_n$ and $x_1,\dotsc,x_n\in\rr$,
\begin{multline}
  \label{eq:airy1fd}
  \pp\!\left(\aipo(t_1)\leq x_1,\dotsc,\aipo(t_n)\leq x_n\right)\\
  =\det\!\left(I-B_0+\bar P_{x_1}e^{(t_2-t_1)\Delta}\bar P_{x_2}e^{(t_3-t_2)\Delta}\dotsm
    \bar P_{x_n}e^{(t_1-t_n)\Delta}B_0\right)_{L^2(\rr)},
\end{multline}
where $B_0$ is given by the kernel $B_0(x,y)=\Ai(x+y)$ defined in \eqref{eq:defB0} and
$\Delta$ is the Laplacian operator. Observe that in all but the last factor of the form
$e^{s\Delta}$ in the above formula it holds that $s>0$, in which case $e^{s\Delta}$ is the
usual heat kernel. This kernel is ill-defined for $s<0$, but it turns out that in this
case the operator $e^{s\Delta}B_0$ makes sense if defined via the integral kernel
\begin{equation}
  e^{s\Delta}B_0(x,y)=e^{2s^3/3+s(x+y)}\Ai(x+y+s^2).\label{eq:sgB0}
\end{equation}
What we mean by this is that if $s,t>0$ then, with this definition the semigroup property
$e^{t\Delta}e^{-s\Delta}B_0=e^{(t-s)\Delta}B_0$ holds.

As we will see in Section \ref{sec:contStat}, it is fruitful to think of the operator
appearing in \eqref{eq:airyfd} as the solution of certain boundary value problem, so we
will refer to formulas like this as \emph{boundary value kernel} formulas. By using
\eqref{eq:sgB0} one can rewrite the definition \eqref{eq:fExtAiry1} of the extended kernel
for $\aipo$ as
\[K^{\rm ext}_1(s,x;t,y)=\begin{dcases*}
  e^{(t-s)\Delta}B_0(x,y) & if $s\geq t$,\\
  -e^{(t-s)\Delta}(I-B_0)(x,y) & if $s<t$.
\end{dcases*}\] It becomes clear then that both the extended kernel formula and the
boundary value kernel formula for Airy$_1$ are obtained from the corresponding formulas
for Airy$_2$ by substituting $H$ with $-\Delta$ and $\K$ with $B_0$. It turns out, as
shown in \cite{bcr}, that the necessary structure behind these formulas hold for a much
wider class of processes, including for instance the stationary GUE Dyson Brownian motion
and non-stationary processes like the Airy$_{2\to1}$ process, and the Pearcey process
\cite{tracyWidom-Pearcey}. For example, for Airy$_{2\to1}$ one has \cite{bcr}
\begin{multline}\label{eq:2to1}
  \pp\!\left(\Bt(t_1)\le x_1,\dotsc,\Bt(t_n)\le x_n\right)\\
  =\det\!\big(I-K_{2\to1}^{t_1}+\bar P_{\tilde x_1}e^{(t_2-t_1)\Delta}\bar P_{\tilde
    x_2}\cdots e^{(t_n-t_{n-1})\Delta}\bar P_{\tilde x_n} e^{(t_1-t_n)\Delta}
  K_{2\to1}^{t_1}\big)_{L^2(\rr)},
\end{multline}
where $\tilde x_i=x_i-t_i^2\uno{t_i\leq0}$ and $K_{2\to1}^{t}(x,y)=K^{\rm ext}_{2\to1}(t,x;t,y)$
with $K^{\rm ext}_{2\to1}$ as in \eqref{eqKCompleteInfinity}.

Interestingly, it is shown in \cite{bcr} that in a setting corresponding to discrete
non-intersecting paths, analogous boundary value kernel formulas can be obtained directly
from applying the Karlin-McGregor formula \eqref{eq:kmcg} (or rather its combinatorial
analogue, the Lindstr\"om-Gessel-Viennot Theorem \cite{lindstrom,gesselViennot}),
bypassing the direct application of the Eynard-Mehta Theorem. In the case of the Airy$_2$
process, a suitable limit of a discrete family of non-intersecting should lead to
\eqref{eq:airyfd} (cf. \eqref{eq:limHrm}). Such a procedure does not seem to work for the
Airy$_1$ process. In fact, in that case the determinantal process used to derive
\eqref{eq:fExtAiry1} is signed (in the sense that the measure defined by the analog of
\eqref{eq:detstr} is signed), see \cite{borFerPrahSasam}, and hence it is not clear how to
associate directly to it a family of non-intersecting paths.

As we will see below, the integral kernels of the operators appearing inside the Fredholm
determinants in \eqref{eq:airyfd}, \eqref{eq:airy1fd} and \eqref{eq:2to1} can be expressed
simply in terms of hitting probabilities of Browian motion. In other words, hitting
probabilities of curves by $\aip$, $\aipo$ and $\Bt$ can be expressed in terms of Fredholm
determinants of the analogous hitting probabilities for Brownian motion. Given the above
discussion (and the discussion in Section \ref{sec:dysonBM}), this is not entirely
surprising in the case of $\aip$, as it follows from the non-intersecting nature of
systems of Brownian paths that can be used to approximate $\aip$. For the same reason, it
is surprising in a sense that the same structure is present in $\aipo$.

\subsection{Continuum statistics and boundary value problems}\label{sec:contStat}

Consider the following problem: compute the probability that, inside a finite interval
$[\ell,r]$, the Airy$_2$ process lies below a given function $g$. The obvious way to
proceed is to take a fine mesh $\ell=t_1<t_2<\dotsm<t_n=r$ of the interval $[\ell,r]$,
take $x_i=g(t_i)$, and attempt to take a limit as $n\to\infty$ in the formula for the
finite dimensional distributions of $\aip$,
\begin{multline}\label{eq:a2mesh}
  \pp\!\left(\aip(t_1)\leq g(t_1),\dotsc,\aip(t_n)\leq g(t_n)\right)\\
  =\det\!\left(I-\K+\bar P_{g(t_1)}e^{(t_1-t_2)H}\bar P_{g(t_2)}e^{(t_2-t_3)H}\dotsm
    \bar P_{g(t_n)}e^{(t_n-t_1)H}\K\right).
\end{multline}
Here the Fredholm determinant is computed on the Hilbert space $L^2(\rr)$, which we will
omit from the subscript in the sequel. By Theorem \ref{thm:Airy2LPP}, $\aip$ has a
continuous version, and hence $\lim_{n\to\infty}\pp\!\left(\aip(t_1)\leq g(t_1),\dotsc,\aip(t_n)\leq
g(t_n)\right)=\pp\!\left(\aip(t)\leq g(t)\text{ for }t\in[\ell,r]\right)$. 
To study the right hand side of \eqref{eq:a2mesh} we need to compute the limit of the
operator appearing inside the determinant. Observe that the last exponential equals
$e^{(r-\ell)H}$, and hence does not depend on $n$. On the other hand, for $s<t$ the
operator $e^{(s-t)H}$ can be thought of as mapping a function $f$ to the solution
$u(t,\cdot)$ at time $t$ of the PDE $\p_tu+Hu=0$ with initial condition
$u(s,\cdot)=f(\cdot)$. Therefore the operator
\begin{equation}
\bar P_{g(t_1)}e^{(t_1-t_2)H}\bar P_{g(t_2)}e^{(t_2-t_3)H}\dotsm e^{(t_{n-1}-t_n)H}
\bar P_{g(t_n)}\label{eq:PbaretH}
\end{equation}
can be thought of as solving the same PDE (backwards in time) on the interval $[\ell,r]$
with the additional condition that all the mass above $g(t_i)$ is removed at each of the
discrete times $t_i$. Note that the PDE is solved backwards because, if we apply this
operator to a function on its right, we first apply $\bar P_{g(t_n)}$, then
$e^{(t_{n-1}-t_n)H}$, then $\bar P_{g(t_{n-1})}$, and so on. Since solving the PDE
$\p_tu+Hu=0$ forward or backwards in time gives the same answer, if we want to think of \eqref{eq:PbaretH}
as being solved forward in time, all we need to do is reverse the order in which the
$g(t_i)$ appear. The result is the following. Given $g\in H^1([\ell,r])$ (i.e. both $g$
and its derivative are in $L^2([\ell,r])$), define an operator $\Theta^g_{[\ell,r]}$
acting on $L^2(\rr)$ as follows: $\Theta^g_{[\ell,r]}f(\cdot)=u(r,\cdot)$, where
$u(r,\cdot)$ is the solution at time $r$ of the boundary value problem
 \begin{equation}
\begin{aligned}
  \p_tu+Hu&=0\quad\text{for }x<g(t), \,\,t\in (\ell,r)\\
  u(\ell,x)&=f(x)\uno{x<g(\ell)}\\
  u(t,x)&=0\quad\text{for }x\ge g(t).
\end{aligned}\label{eq:bdval}
\end{equation}
Further, define $\hat g(t)=g(\ell+r-t)$. Then
\begin{equation}
  \left\|\bar P_{g(t_1)}e^{(t_1-t_2)H}\bar P_{g(t_2)}e^{(t_2-t_3)H}\dotsm
    \bar P_{g(t_n)}e^{(t_n-t_1)H}\K-\Theta^{\hat g}_{[\ell,r]}e^{(t_n-t_1)H}\K\right\|_1
  \xrightarrow[n\to\infty]{}0.\label{eq:discrToCont}
\end{equation}
Since the convergence holds in trace class norm,  \eqref{eq:discrToCont} can be used to answer the question
with which we started this subsection:

\begin{thm}[\cite{cqr}, Theorem \grefn{thm:aiL}]\label{thm:aiL}
  \begin{equation}\pp\!\left(\aip(t)\leq g(t)\text{ for
      }t\in[\ell,r]\right)=\det\!\left(I-K_{\Ai}+\Theta^g_{[\ell,r]}  e^{(r-\ell)H}K_{\Ai}\right).\label{eq:basic}
  \end{equation}
\end{thm}

Observe that we have written $g$ instead of $\hat g$ in \eqref{eq:basic}. We may do this
because the Airy$_2$ is invariant under time reversal, so we can replace $g$ by $\hat g$
on the left hand side.

The limit \eqref{eq:discrToCont} is proved in Proposition \gref{prop:theta} (in fact only
along the dyadic sequence $n_k=2^k$, but this is enough for deducing Theorem \ref{thm:aiL}). The
proof is based on the following probabilistic representation of the solutions of the
above boundary value problem: if
$\Theta^g_{[\ell,r]}(x,y)$ denotes the integral kernel of $\Theta^g_{[\ell,r]}$, then
\begin{multline}\label{eq:ThetaL}
  \Theta^g_{[\ell,r]}(x,y)=e^{\ell x-ry+(r^3-\ell^3)/3}\frac{e^{-(x-y)^2/4(r-\ell)}}{\sqrt{4\pi(r-\ell)}}\\
  \cdot\pp_{\hat b(\ell)=x-\ell^2,\hat b(r)=y-r^2}\!\left(\hat b(s)\leq g(s)-s^2\text{ on
    }[\ell,r]\right),
\end{multline}
where the probability is computed with respect to a Brownian bridge $\hat b(s)$ from
$x-\ell^2$ at time $\ell$ to $y-r^2$ at time $r$ and with diffusion coefficient $2$. This
formula is Theorem \gref{thm:thetaLgen}, its proof is based on an application of the
Feynman-Kac and Cameron-Martin-Girsanov formulas.
 
 The argument that proves Theorem \ref{thm:aiL} can be adapted to obtain a similar result
 for Airy$_1$. Fix $\ell<r$. Given $g\in H^1([\ell,r])$, define an operator
 $\Lambda^g_{[\ell,r]}$ acting on $L^2(\rr)$ as follows:
 $\Lambda^g_{[\ell,r]}f(\cdot)=u(r,\cdot)$, where $u(r,\cdot)$ is the solution at time $r$
 of the boundary value problem
 \begin{equation}
\begin{aligned}
  \p_tu-\Delta u&=0\quad\text{for }x<g(t), \,\,t\in (\ell,r)\\
  u(\ell,x)&=f(x)\uno{x<g(\ell)}\\
  u(t,x)&=0\quad\text{for }x\ge g(t).
\end{aligned}\label{eq:bdvalo}
\end{equation}

\begin{thm}[\cite{quastelRemAiry1}, Theorem \orefn{thm:aiLo}]\label{thm:aiLo}
  \begin{equation}\label{eq:basic1}
    \pp\!\left(\aipo(t)\leq g(t)\text{ for
      }t\in[\ell,r]\right)=\det\!\left(I-B_0+\Lambda^g_{[\ell,r]}e^{-(r-\ell)\Delta}B_0\right).
  \end{equation}
\end{thm}

Although the proof of this result is similar to the proof for the Airy$_2$ case, the
argument is a bit more involved because, as written, the operator in the determinant is
not trace class, so one needs to conjugate appropriately. 
Of course, similar arguments should allow one to obtain continuum statistics formulas for
other processes for which boundary value kernel formulas are available (see
\cite{bcr} for the case of stationary GUE Dyson Brownian motion).

The operator
$\Lambda^g_{[\ell,r]}$ also has a simple representation in terms of Brownian motion (see
\cite{quastelRemAiry1}), which has recently been used in \cite{ferrariFrings} to verify
numerically the experimental values obtained in \cite{takeuchiSano2} for the persistence
probabilities of Airy$_1$.  The negative persistence exponent is defined by
\begin{equation}
\pp(\aipo(t)\le m, ~0\le t\le L) \sim e^{-\kappa_- L}
\end{equation}
where $m$ is the mean of $F_{\rm GOE}$.  
Takeuchi has measured $\kappa_-\approx 3.2\pm 0.2$ in computer simulations of 
the Eden model \cite{takeuchiSano2}.  
Ferrari and Frings \cite{ferrariFrings}  have computed numerically
 \eqref{eq:basic1} finding
\begin{equation}\kappa\approx 2.9,\end{equation}
which is fairly close.  Note that 
Takeuchi has also measured the positive persistence probabilities 
$\pp(\aipo(t)\ge m, ~0\le t\le L) \sim e^{-\kappa_+ L}$
An interesting question  is whether there exists a simple enough mathematical formula to
check such a thing.

\section{Applications}
\label{sec:appl}

In this section we will describe some applications of the boundary value kernel formulas
for Airy processes which were introduced in the previous section. The first two
applications refer to asymptotic statistics for directed polymers and LPP, while the next
two involve respectively the Airy$_1$ and Airy$_{2\to1}$ processes.

\subsection{Point-to-line LPP and GOE}
\label{sec:goe}

Recall the variational formula \eqref{eq:Airy2GOE} relating the Airy$_2$ process with the
Tracy-Widom GOE distribution:
\begin{equation}
\pp\!\left(\sup_{x\in\rr}\big\{\aip(x)-x^2\big\}\leq m\right)=F_{\rm GOE}(4^{1/3}m).\label{eq:johGOE}
\end{equation}
As we explained in Section \ref{sec:lpp}, Johansson's proof \cite{johansson} was very
indirect, relying on the convergence of the spatial fluctuations of point-to-point LPP to
$\aip$ together with \eqref{eq:lpplptline} and \eqref{eq:LPPGOE}.

A direct proof of this variational formula was provided in \cite{cqr}, based on the
continuum statistics formula given in Theorem \ref{thm:aiL}. An interesting consequence of
this derivation was that it allowed to identify the factor of $4^{1/3}$ on the right hand
side of the identity, which had been lost in Johansson's argument in the process of
translating between the available results at the time (see Section \gref{sec:indirect} for
an account of how to get the correct factor directly from LPP).

We will explain next the derivation of the formula, skipping some details. We rewrite
the desired probability as
\[\lim_{L\to\infty}\pp\!\left(\aip(t)\leq m+t^2~~\forall\,t\in[-L,L]\right).\]
For fixed $L>0$, Theorem \ref{thm:aiL} implies that this probability is given by
\begin{equation}
\det\!\left(I-K_{\Ai}+\Theta_Le^{2LH}K_{\Ai}\right),\label{eq:detAiry2}
\end{equation}
where
\[\Theta_L=\Theta^{g(t)=t^2+m}_{[-L,L]}.\]
The nice thing is that the choice of $g(t)=t^2+m$ is the simplest possible from the point
of view of explicit calculations, because it cancels exactly the parabola appearing on the
right hand side of \eqref{eq:ThetaL}. The probability appearing in that formula is then
reduced to the probability of a Brownian bridge staying below level $m$, and this is easy
to compute using the reflection principle (method of images):
\begin{equation}
\begin{split}
  &\pp_{\hat b(-L)=x-L^2,\hat b(L)=y-L^2}\!\left(\hat b(s)\leq m\text{ on }[-L,L]\right)\\
  &\hspace{1in}=1-\pp_{\hat b(-L)=x-L^2,\hat b(L)=y-L^2}\!\left(\hat b(s)>m\text{ for some }s\in[-L,L]\right)\\
  &\hspace{1in}=1-e^{-(x-m-L^2)(y-m-L^2)/2L}
\end{split}\label{eq:refl}
\end{equation}
(we leave the simple computation to the reader, alternatively see page 67 of
\cite{handbookBM}). Putting this back in $\Theta_L$ gives
\begin{equation}
\Theta_L=\bar P_{m+L^2}e^{-2LH}\bar P_{m+L^2}-\bar P_{m+L^2}R_L\bar P_{m+L^2},\label{eq:thetaL}
\end{equation}
where $R_L$ is the reflection term
\begin{equation}\label{eq:RL}
  R_L(x,y)=\frac{1}{\sqrt{8\pi L}}e^{-(x+y-2m-2L^2)^2/8L-(x+y)L+2L^3/3}.
\end{equation}
The $e^{-2LH}$ in the first term in $\Theta_L$ comes from the 1 in \eqref{eq:refl} and
appears by either reversing the use of the Cameron-Martin-Girsanov and Feynman-Kac
formulas in the derivation of \eqref{eq:thetaL} or by an explicit computation of the
integral kernel of $e^{-(r-\ell)H}$ as
\[e^{-(r-\ell)H}(x,y)=e^{\ell x-ry+(r^3-\ell^3)/3}\frac{e^{-(x-y)^2/4(r-\ell)}}{\sqrt{4\pi(r-\ell)}}.\]

Referring to \eqref{eq:detAiry2}, we have by the cyclic property of determinants \eqref{eq:cyclic} and the identity
$e^{2LH}K_{\Ai}=(e^{LH}K_{\Ai})^2$ (which follows from Remark \ref{airyrem}) that
\begin{equation}\label{eq:basiccyclic}
  \pp\!\left(\aip(t)\leq t^2+m\text{ for  }t\in[-L,L]\right)=\det\!\left(I-K_{\Ai}+e^{LH}K_{\Ai}
    \Theta_L e^{LH}K_{\Ai}\right).
\end{equation}
To obtain the $L\to \infty$ asymptotics, we decompose $\Theta_L$ so as to expose the two
limiting terms, as well as a remainder term $\Omega_L$:
\begin{equation}
\Theta_L=e^{-2LH}-R_L+\Omega_L,\label{eq:OmegaL}
\end{equation}
where $\Omega_L=\big(R_L-\bar P_{m+L^2}R_L\bar P_{m+L^2}\big)
-\big(e^{-2LH}-\bar P_{m+L^2}e^{-2LH}\bar P_{m+L^2}\big)$.  It is shown in \cite{cqr}
that
\begin{equation}
\left\|e^{LH}K_{\Ai}\Omega_Le^{LH}K_{\Ai}\right\|_1\xrightarrow[L\to\infty]{}0.\label{eq:OmegaLto0}
\end{equation}
The proof amounts essentially to asymptotic analysis involving the Airy function. In view
of this fact and the decomposition \eqref{eq:OmegaL}, and since
$e^{LH}K_{\Ai}e^{-2LH}e^{LH}K_{\Ai}=K_{\Ai}$, we see that the key point is the limiting
behaviour in $L$ of
\[e^{LH}K_{\Ai}R_Le^{LH}K_{\Ai}.\]

To explain how this last product can be computed we will proceed in a slightly formal
manner through an argument based on the Baker-Campbell-Hausdorff formula, as done for a
related problem in \cite{qr-airy1to2} (see Section \ref{sec:airy2to1}). Since $\K$ is a
projection and $H$ leaves $\K$ invariant, we will pretend that $e^{LH}$ and $\K$
commute, so we have to compute the limit of $e^{LH}R_{L}e^{LH}$.  Define the
\emph{reflection operator} $\varrho_m$ by
\[\varrho_mf(x)=f(2m-x).\]
Then the operator $R_{L}$ defined in \eqref{eq:RL} can be rewritten as
\begin{equation}
R_{L}=e^{(2L^3)/3}e^{-L\xi}\varrho_{m+L^2}e^{2L\Delta}e^{-L\xi}
=e^{(2L^3)/3}e^{-L\xi}e^{L\Delta}\varrho_{m+L^2}e^{L\Delta}e^{-L\xi}.\label{eq:Rgauss}
\end{equation}
Here $e^{r\xi}$ ($\xi$ stands for a generic variable) denotes the multiplication operator
$(e^{r\xi}f)(x)=e^{rx}f(x)$. The second equality follows from the reflection principle
applied to the heat kernel. 

The following identities will be useful, where $[\cdot,\cdot]$ denotes commutator:
\[[H,\Delta]=[\xi,\Delta]=-2\nabla,\qquad[H,\nabla]=[\xi,\nabla]=-I,\qquad[H,\xi]=-2\nabla.\]
If $A$ and $B$ are two operators such that $[A,[A,B]]=c_1I$ and $[B,[A,B]]=c_2I$ for some
$c_1,c_2\in\rr$, then the Baker-Campbell-Hausdorff formula reads
\begin{equation}
  \label{eq:bch}
  e^{A}e^{B}=e^{A+B+\frac12[A,B]+\frac1{12}[A,[A,B]]-\frac1{12}[B,[A,B]]}.
\end{equation}
Using this we have
\[e^{-L\xi}e^{L\Delta}=e^{L^3/6}e^{L\Delta+L^2\nabla-L\xi}.\]
Using  the Baker-Campbell-Hausdorff formula again we deduce that
\[e^{LH}e^{-L\xi}e^{L\Delta}=e^{L^3/6}e^{LH}e^{L\Delta+L^2\nabla-L\xi}
  =e^{-L^3/3}e^{L^2\nabla},\]
while an analogous computation yields
\[e^{L\Delta}e^{-L\xi}e^{LH}=e^{-L^3/3}e^{-L^2\nabla}.\] Employing these identities on the
right hand side of \eqref{eq:Rgauss} yields
\[e^{LH}R_{L}e^{LH}=e^{L^2\nabla}\varrho_{m+L^2}e^{-L^2\nabla}.\]
Since $e^{r\nabla}$ is the shift operator $(e^{r\nabla}f)(x)=f(x+r)$, we have
$e^{r\nabla}\varrho_m=\varrho_me^{-r\nabla}=\varrho_{m-r/2}$, and we obtain
\[e^{LH}R_{L}e^{LH}=\varrho_{m}.\]
Remarkably, the result does not depend on $L$. The conclusion from using this,
\eqref{eq:OmegaLto0} and \eqref{eq:basiccyclic} in \eqref{eq:OmegaL} and taking
$L\to\infty$ is that
\begin{equation}
  \label{eq:resultLGOE}
   \pp\!\left(\aip(t)\leq t^2+m\text{ for all  }t\in\rr\right)=\det\!\left(I-K_{\Ai}\varrho_m\K\right).
\end{equation}
The use of the Baker-Campbell-Formula in the derivation of this identity can be replaced by an
explicit integral calculation (see the proof of Proposition \gref{thm:goe}).

To finish our proof of \eqref{eq:johGOE} we need to show that the right hand side of
\eqref{eq:resultLGOE} equals $F_{\rm GOE}(4^{1/3}m)$. Recall the definition of the kernel
$B_0(x,y)=\Ai(x+y)$ and observe that $\K=B_0P_0B_0$. Recall also that the shifted Airy
functions form a generalized orthonormal basis of $L^2(\rr)$ (see Remark \ref{airyrem}),
which implies that $B_0^2=I$. Therefore we can use the cyclic property of determinants \eqref{eq:cyclic} to deduce that
\[\det\!\left(I-K_{\Ai}\varrho_m\K\right)=\det\!\left(I-P_0B_0\varrho_mB_0P_0\right).\]
Now
\[B_0\varrho_mB_0(x,y)=\int_{-\infty}^\infty d\lambda\Ai(x+\lambda)\Ai(2m-\lambda+y),\]
and using the identity \eqref{eq:airyIdentity} we deduce that
\begin{equation}
  \label{eq:wtBm}
  B_0\varrho_mB_0(x,y)=\widetilde B_m(x,y):=2^{-1/3}\Ai(2^{-1/3}(x+y+2m),
\end{equation}
and thus
\[\det\!\left(I-K_{\Ai}\varrho_m\K\right)=\det\!\left(I-P_0\widetilde B_mP_0\right).\]
Performing the change of variables $x\mapsto2^{1/3}x$, $y\mapsto2^{1/3}y$ in the series
defining the last Fredholm determinant shows that the determinant on the right hand side
of \eqref{eq:resultLGOE} equals $\det(I-P_0B_{4^{1/3}m}P_0)$, which is $F_{\rm
  GOE}(4^{1/3}m)$ by \eqref{eq:GOE}.

\subsection{Endpoint distribution of directed polymers}
\label{sec:endpoint}

In the setting of geometric LPP (see Section \ref{sec:lpp}), consider the random variable
\[\kappa_N=\min\!\left\{k\in\{-N,\dotsc,N\}\!:\sup_{j=-N,\dotsc,k}L^{\rm
    point}_N(j)=\sup_{j=-N,\dotsc,N}L^{\rm point}_N(j)\right\}.\] $\kappa_N$ corresponds
to the location of the endpoint of the maximizing path in point-to-line LPP.

Interest in the scaling properties and distribution of this random variable goes back at
least to the early 1990's. One can also consider the analogous random variable in the
setting of directed random polymers, but due to the KPZ universality conjecture one
expects that the asymptotic behavior and statistics are the same as in
LPP. \citet{mezardParisi} considered the polymer case and derived non-rigorously the
scaling relation
\begin{equation}
  \label{eq:endpointScaling}
  |\kappa_N|\sim N^{2/3}
\end{equation}
(c.f. \eqref{eq:ZNPE}). In view of this we define the rescaled endpoint
\begin{equation}
  \label{eq:TN}
  \ct_N=c_3^{-1}N^{-2/3}\kappa_N,
\end{equation}
where $c_3$ is the constant appearing in \eqref{eq:rescaledptp}. Recalling the definition
of the rescaled point-to-point last passage time \eqref{eq:rescaledptp} as the linear
interpolation of the values given by
\[H^{\rm point}_N(t)=\frac1{c_2N^{1/3}}\left[L^{\rm
    point}(N+c_3^{-1}N^{-2/3}t,N-c_3^{-1}N^{-2/3}t)-c_1N\right]\] for $t$ such that
$c_3^{-1}N^{-2/3}t\in\{-N,\dotsc,N\}$ we deduce that
\[\ct_N=\min\!\left\{t\in\rr\!:\sup_{s\leq t}H^{\rm
    point}_N(s)=\sup_{s\in\rr}H^{\rm point}_N(s)\right\}.\] Recalling that $H_N(t)$
converges to $\aip(t)-t^2$ by Theorem \ref{thm:Airy2LPP} it becomes clear that $\ct_N$
should converge to the point where $\aip(t)-t^2$ attains it maximum. In fact, this is what
Johansson proved, although he had to make a (very reasonable) technical assumption on the
Airy$_2$ process which he was not able to prove with the tools available at the time:

\begin{thm}[\cite{johansson}]\label{thm:TNcvgce}
  Assume that the process $\aip(t)-t^2$ attains its maximun at a unique point and let
  \[\ct=\argmax_{t\in\rr}\big\{\aip(t)-t^2\big\}.\]
  Then
  \[\ct_N\xrightarrow[N\to\infty]{}\ct\]
  in the sense of convergence in distribution.
\end{thm}

Although the result is of course very interesting, as it shows that the limiting endpoint
distribution exists (under the technical assumption), it gives no information on the
distribution of $\ct$. Quoting \citet{johanssonRMandDetPr}, for all we know $\ct$ could be
Gaussian. Nevertheless, from KPZ universality one expects that this is not the case. For
example, \citet{halpZhang} conjectured on the basis of analogy with the argmax of a
Brownian motion minus a parabola (for which one has a complete analytical solution, see
\cite{groeneboom}), that the tails of $\ct$ decay like $e^{-ct^3}$, which of course rules
out Gaussian behavior.

It turns out that the distribution of $\ct$ can be computed explicitly through an argument
based on the continuum statistics formula of Theorem \ref{thm:aiL}. This was done in
\cite{mqr}, where in fact the joint density of
\[\ct=\argmax_{t\in\rr}\big\{\aip(t)-t^2\big\}\qqand\cm=\max_{t\in\rr}\big\{\aip(t)-t^2\big\}\]
was computed. Moreover, the argument implies that the maximum of $\aip(t)-t^2$ is attained
at a unique point, thus completing the proof of Theorem \ref{thm:TNcvgce}. The uniqueness
of the maximum was also proved slightly earlier by \citet{corwinHammondKPZ} using
completely different techniques, and a proof for general stationary processes is now
available \cite{pimentelStatMax}.

The computation is as follows. For simplicity we will assume the uniqueness of the
maximizing point of $\aip(t)-t^2$, and will explain later how the uniqueness can actually
be obtained from this argument. Let $(\cm_L,\ct_L)$ denote the maximum and the location of
the maximum of $\aip(t)-t^2$ restricted to $t\in[-L,L]$, and let $f_L$ be the joint
density of $(\cm_L,\ct_L)$.  By results of \cite{corwinHammond}, the joint density
$f(m,t)$ of $\cm, \ct$ is well approximated by $f_L(m,t)$,
\[f(t,m) =\lim_{L\to \infty} f_L(t,m).\] By definition,
\[f_L(t,m) = \lim_{\delta\to0}\lim_{\ep\to0}\frac1{\ep\delta}\pp\!\left(\cm_L \in
  [m,m+\ep],\, \ct_L\in [t,t+\delta]\right),\] provided that the limit exists. The main
contribution in the above expression comes from paths entering the space-time box
$[t,t+\delta]\times[m,m+\ep]$ and staying below the level $m$ outside the time interval
$[t,t+\delta]$. More precisely, if we denote by $\underline D_{\ep,\delta}$ and $\overline
D_{\ep,\delta}$ the sets
\begin{equation}
  \begin{aligned}
    \underline{D}_{\ep,\delta} &= \Big\{\aip(s)-s^2\leq m,\,\,\,
    s\in[t,t+\delta]^\text{c},\,\aip(s)-s^2\leq m+\ep,\,\,\,s\in[t,t+\delta], \\
    &\hspace{2.6in}\aip(s)-s^2 \in [m,m+\ep]\text{ for some }s\in [t,t+\delta]\Big\},\\
    \shortintertext{and} \overline{D}_{\ep,\delta}& = \Big\{\aip(s)-s^2\leq m+\ep,\,\,\,
    s\in[-L,L],\,\aip(s)-s^2 \in [m,m+\ep]\text{ for some }s\in [t,t+\delta]\Big\},
  \end{aligned}
\end{equation}
then
\begin{equation}
  \underline D_{\ep,\delta}\subseteq\left\{\cm_L \in
    [m,m+\ep],\, \ct_L\in [t,t+\delta]\right\}\subseteq\overline D_{\ep,\delta}.
\end{equation}
Letting $\underline
f(t,m)=\lim_{\delta\to0}\lim_{\ep\to0}\frac1{\ep\delta}\pp\big(\underline
D_{\ep,\delta}\big)$ and defining $\overline f(t,m)$ analogously (with $\overline
D_{\ep,\delta}$ instead of $\underline D_{\ep,\delta}$) we deduce that $\underline
f(t,m)\leq f(t,m)\leq\overline f(t,m)$. In what follows we will compute $\underline
f(t,m)$. It will be clear from the argument that for $\overline f(t,m)$ we get the same
limit. The conclusion is that
\[f_L(t,m)=\lim_{\delta\to0}\lim_{\ep\to0}\frac1{\ep\delta}\pp\big(\underline
D_{\ep,\delta}\big).\]

We rewrite this last equation as
\begin{equation}\label{first}
  f_L(t,m)= \lim_{\delta\to0}\lim_{\ep\to0}\frac1{\ep\delta}\Big[\,\pp\!\left(\aip(s)\leq
    h_{\ep,\delta}(s),\,s\in[-L,L]\right)-\pp\!\left(\aip(s)\leq
    h_{0,\delta}(s),\,s\in[-L,L]\right)\Big],
\end{equation}
where
\[h_{\ep,\delta}(s)=s^2+m+\ep\uno{s\in[t,t+\delta]}.\] These two probabilities have
explicit Fredholm determinant formulas by Theorem \ref{thm:aiL}. We get, using the cyclic
property of determinants as in \eqref{eq:basiccyclic},
\begin{multline}
  f_L(t,m) = \lim_{\delta\to0}\lim_{\ep\to0}\frac1{\ep\delta} \left[
    \det\!\left(I-K_{\Ai}+e^{LH}K_{\Ai}\Theta^{h_{\ep,\delta}}_{[-L,L]}e^{LH}K_{\Ai}\right) \right.\\
  \left.-\det\!\left(I-K_{\Ai}+e^{LH}K_{\Ai}\Theta^{h_{0,\delta}}_{[-L,L]}e^{LH}K_{\Ai}\right)\right].
\end{multline}
The limit in $\ep$ becomes a derivative
\begin{align}
  f_L(t,m) &=\lim_{\delta\to0}\frac{1}{\delta}\,\p_\beta
  \!\left.\det\!\left(I-K_{\Ai}+e^{LH}K_{\Ai}\Theta^{h_{\beta,\delta}}_{[-L,L]}e^{LH}K_{\Ai}\right)\right|_{\beta=0},
\end{align}
which in turn gives a trace by \eqref{derivativeofadet},
\begin{multline}\label{grd}
  f_L(t,m)=\det\!\left(I-K_{\Ai}+e^{LH}K_{\Ai}\Theta^{h_{0,\delta}}_{[-L,L]}e^{LH}K_{\Ai}\right)\\
  \cdot\lim_{\delta\to0}\frac{1}{\delta}
  \tr\!\left[(I-K_{\Ai}+e^{LH}K_{\Ai}\Theta^{h_{0,\delta}}_{[-L,L]}e^{LH}K_{\Ai})^{-1}
    e^{LH}K_{\Ai}\left[\p_\beta\Theta^{h_{\beta,\delta}}_{[-L,L]}\right]_{\beta=0}e^{LH}K_{\Ai}\right].
\end{multline}
One has to check here that the required limits hold in trace class norm, see
\cite{mqr}. Note that $h_{0,\delta}=g_m$, where $g_m$ is the parabolic barrier
\[g_m(s)=s^2+m,\] so in particular the determinant and the first factor inside the trace
do not depend on $\delta$. We know moreover from the arguments in Section \ref{sec:goe}
that
\begin{equation} \label{eq:GOEcvg} \lim_{L\to\infty}\left(
    I-K_{\Ai}+e^{LH}K_{\Ai}\Theta^{h_{0,\delta}}_{[-L,L]}e^{LH}K_{\Ai}\right)=I-\K\varrho_m\K
\end{equation}
in trace norm. In particular, we have
\[\lim_{L\to\infty}
\det\!\left(I-K_{\Ai}+e^{LH}K_{\Ai}\Theta^{h_{0,\delta}}_{[-L,L]}e^{LH}K_{\Ai}\right) =
F_\mathrm{GOE}(4^{1/3}m).\]

The next step is to compute $\p_\beta\Theta^{h_{\beta,\delta}}_{[-L,L]}\!\mid_{\beta=0}$.
Recalling that $h_{0,\delta}(s)=g_m(s)= s^2+m$ and also $h_{\ep,\delta}(s)=g_{m+\ep}(s)$
for $s\in[t,t+\delta]$ we have, by the semigroup property,
\[\Theta^{h_{\ep,\delta}}_{[-L,L]}-\Theta^{h_{0,\delta}}_{[-L,L]}
=\Theta^{g_m}_{[-L,t]}\left[\Theta^{g_{m+\ep}}_{[t,t+\delta]}
  -\Theta^{g_m}_{[t,t+\delta]}\right]\Theta^{g_m}_{[t+\delta,L]}.\] Computing the desired
derivative involves just the middle bracket, which we note corresponds to the same
boundary value problem as in Section \ref{sec:goe}, only at two different levels $m$ and
$m+\ep$. Since we have explicit formulas, the derivative can be computed explicitly. The
computation is slightly tedious, and the only delicate part is to justify that the
necessary limits occur in trace class norm, we refer to \cite{mqr} for the details.

Going back to \eqref{grd}, we recall that the trace is linear and continuous under the
trace class norm topology, so in view of the preceding discussion we have
\begin{multline}
  \lim_{L\to\infty}f_L(t,m)\\=F_{\rm
    GOE}(4^{1/3}m)\tr\!\left[(I-K_{\Ai}\varrho_mK_{\Ai})^{-1}\lim_{L\to\infty}\lim_{\delta\to0}\frac{1}{\delta}
    e^{LH}K_{\Ai}\left[\p_\beta\Theta^{h_{\beta,\delta}}_{[-L,L]}\right]_{\beta=0}e^{LH}K_{\Ai}\right].\label{eq:grd2}
\end{multline}
Once again we need to compute limits, again taking care that they hold in trace class norm
as necessary. We skip the details and just write down the result, \peqref{eq:grt}:
\begin{equation}\label{eq:grt}
  \lim_{L\to\infty}\lim_{\delta\to
    0}\frac{1}{\delta}e^{LH}K_{\Ai}\left[\p_\beta\Theta^{h_{\beta,\delta}}_{[-L,L]}\right]_{\beta=0}
  e^{LH}K_{\Ai}=\Psi,
\end{equation}
where
\[\Psi(x,y)=B_0P_0\psi_{t,m}(x)B_0P_0\psi_{-t,m}(y)\]
and
\[\psi_{t,m}(x)=2e^{t^3+(m+x)t}\bigg[\Ai'(m+t^2+x)+t\Ai(m+t^2+x)\bigg]\]
(we remark that we have written these formulas in a slightly different way compared to
\cite{mqr}, but the reader should have no problem translating between the formulas). The
limit in $\delta$ is relatively straightforward, while the limit in $L$ involves an
argument similar to the one used in Section \ref{sec:goe}. Using this formula in
\eqref{eq:grd2} the trace becomes
\[\tr\!\left[ (I-\K\varrho_m\K)^{-1}\Psi\right] =\left\langle
  (I-\K\varrho_m\K)^{-1}B_0P_0\psi_{t,m},B_0P_0\psi_{-t,m}\right\rangle,\] where
$\langle\cdot,\cdot\rangle$ denotes inner product in $L^2(\rr)$.
    
It only remains to simplify the expression. We first use \eqref{eq:wtBm} and the facts
that $\K=B_0P_0B_0$, $B_0^2=I$ and $B_0^*=B_0$ to write
\begin{align}
  \left\langle(I-\K\varrho_m\K)^{-1}B_0P_0\psi_{t,m},B_0P_0\psi_{-t,m}\right\rangle &=
  \left\langle(I-B_0P_0
    B_mP_0B_0)^{-1}B_0P_0\psi_{t,m},B_0P_0\psi_{-t,m}\right\rangle\\
  &= \left\langle B_0(I-P_0
    B_mP_0)^{-1}P_0\psi_{t,m},B_0P_0\psi_{-t,m}\right\rangle\\
  &= \left\langle(I-P_0 B_mP_0)^{-1}P_0\psi_{t,m},P_0\psi_{-t,m}\right\rangle.
\end{align}
Next we introduce the scaling operator $Sf(x)=f(2^{1/3}x)$. One can check easily that
$S^{-1}=2^{1/3}S^*$ and that $P_0$ commutes with $S$ and $S^{-1}$. We also have $S
B_mS^{-1}=B_{4^{1/3}m}$.  Thus writing $\tilde m=2^{-1/3}m$ we get
\begin{align}
  \Big\langle(I-P_0 B_mP_0)^{-1}P_0\psi_{t,m},P_0\psi_{-t,m}\Big\rangle
  &=\left\langle(I-S^{-1}P_0B_{2\tilde m}P_0S)^{-1}P_0\psi_{t,m},P_0\psi_{-t,m}\right\rangle\\
  &=\left\langle S^{-1}(I-P_0B_{2\tilde m}P_0)^{-1}P_0S\psi_{t,m},P_0\psi_{-t,m}\right\rangle\\
  &=2^{1/3}\left\langle(I-P_0B_{2\tilde
      m}P_0)^{-1}P_0S\psi_{t,m},P_0S\psi_{-t,m}\right\rangle.
\end{align} which is equal to $2^{1/3}\gamma(t,4^{1/3}m)$.

Using this formula in \eqref{eq:grd2} yields the joint density of $\ct$ and $\cm$.  Define
the resolvent kernel
\[\varsigma_m(x,y)=(I-P_0B_mP_0)^{-1}(x,y)\]
and, for $t,m\in\rr$, define
\[\Psi_{t,m}(x,y)=2^{1/3}\psi_{t,m}(2^{1/3}x)\psi_{-t,m}(2^{1/3}y)\]
and
\begin{equation}
  \label{eq:gamma}
  \gamma(t,m)=2^{1/3}\int_0^\infty dx\int_0^\infty dy\,\psi_{-t,4^{-1/3}m}(2^{1/3}x)\varsigma_{m}(x,y)\psi_{t,4^{-1/3}m}(2^{1/3}y).
\end{equation}
 
\begin{thm}[\cite{mqr}, Theorem \prefn{thm:fPE}]\label{thm:fPE}
  The joint density $f(t,m)$ of $\ct$ and $\cm$ is given by
  \begin{equation}\label{eq:eff}
    \begin{split}
      f(t,m)&=\gamma(t,4^{1/3}m)F_\mathrm{GOE}(4^{1/3}m)\\
      &=\det\!\big(I-P_0B_{4^{1/3}m}P_0+P_0\Psi_{t,m}P_0\big)-F_\mathrm{GOE}(4^{1/3}m).
    \end{split}
  \end{equation}
\end{thm}

To see where the second equality in \eqref{eq:eff} comes from, observe that
$\gamma(t,4^{1/3}m)$ equals the trace of the operator
$(I-P_0B_{4^{1/3}m}P_0)^{-1}P_0\Psi_{t,m}P_0$ and that $\Psi_{t,m}$ is a rank one
operator. The identity now follows that from the general fact that for two operators $A$
and $B$ such that $B$ is rank one, one has
$\det(I-A+B)=\det(I-A)\big[1+\tr\big((I-A)^{-1}B\big)\big]$.

\begin{figure}
  \centering
  \includegraphics[width=4.5in]{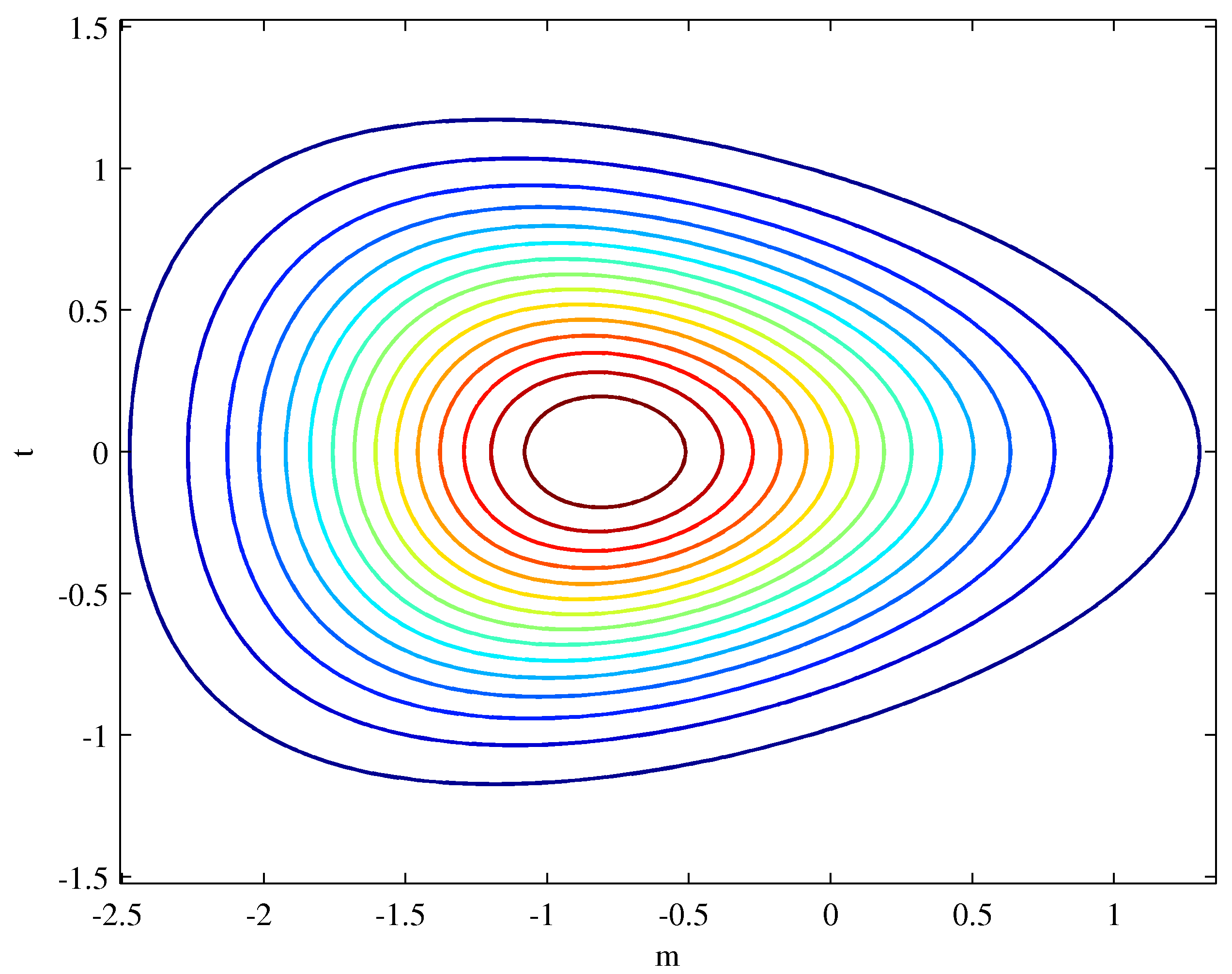}
  \vspace{-14pt}
  \caption{Contour plot of the joint density of $\cm$ and $\ct$.}
  \label{fig:contour}
\end{figure}

\begin{figure}
  \centering \hspace{-0.0in}\includegraphics[width=6in]{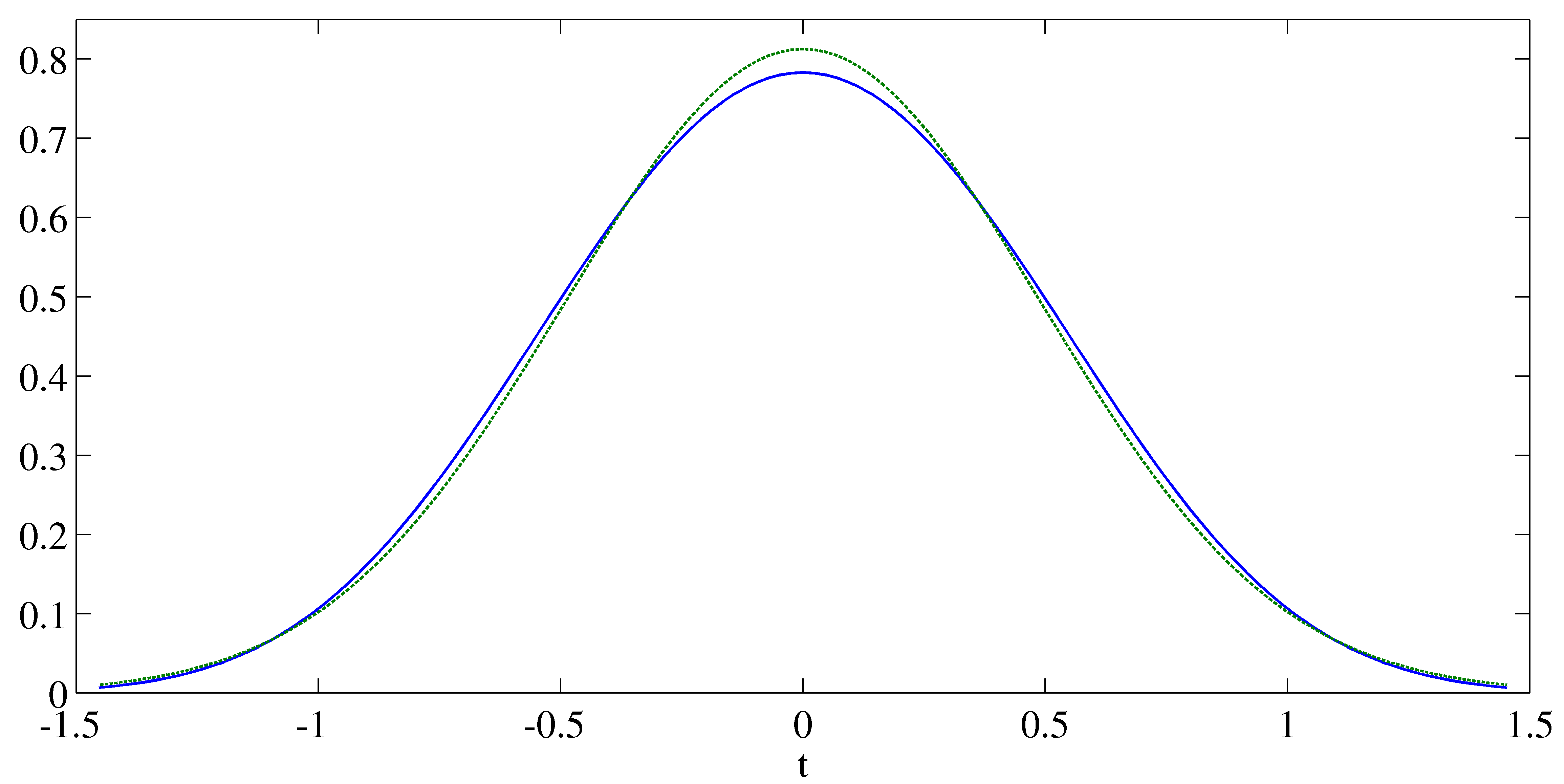} \vspace{-14pt}
  \caption{Plot of the density of $\ct$ compared with a Gaussian density with the same
    variance 0.2409 (dashed line). The excess kurtosis $\ee(\ct^4)/\ee(\ct^2)^2-3$ is
    $-0.2374$.}
  \label{fig:density}
\end{figure}

Integrating over $m$ one obtains a formula for the probability density $f_{\rm end}(t)$ of
$\ct$. Unfortunately, it does not appear that the resulting integral can be calculated
explicitly, so the best formula one has is
\begin{equation}
  f_{\rm end}(t)= \int_{-\infty}^\infty dm\,f(t,m).
\end{equation}
One can readily check nevertheless that $f_{\rm end}(t)$ is symmetric in $t$. The second
formula for $f(t,m)$ is suitable for numerical computations, using the numerical scheme
and Matlab toolbox developed by Bornemann in \cite{bornemann2,bornemann1} for the
computation of Fredholm determinants. Figure \ref{fig:contour} shows a contour plot of the
joint density of $\cm$ and $\ct$, while Figure \ref{fig:density} shows a plot of the
marginal $\ct$ density.

As we mentioned, interest in this problem dates back at least two decades. In particular,
there has been a resurgence of interest in the last couple of years. An alternative way to
obtain the Airy${}_2$ process is as a limit in large $N$ of the top path in a system of
$N$ non-intersecting random walks, or Brownian motions, the so called vicious walkers
\cite{fisher} (this is of course related to the setting presented in Section
\ref{sec:dysonBM}). \cite{SMCR,feierl2,RS1,RS2} obtained various expressions for the
joint distributions of $\cm$ and $\ct$ in such a system at finite $N$.  \cite{forrester}
obtained the $F_{\rm GOE}$ distribution from large $N$ asymptotics non-rigorously, and
furthermore made connections between these problems and Yang-Mills theory. But for several
years people were not able to perform asymptotic analysis on the formulas obtained for
$\ct$ at finite $N$.

After \cite{mqr} appeared, \citet{schehr} succeeded in extracting
asymptotics from the vicious walkers formula, and obtained an alternative formula for
$f(t,m)$. His formula is given as follows. The Painlev\'{e} II equation
\eqref{in:10}-\eqref{in:11} has a Lax pair formulation
\begin{equation}
\frac{\partial}{\partial \zeta} \Phi = A \Phi, \qquad \; \frac{\partial}{\partial s} \Phi = B \Phi\label{in:8}
\end{equation}
for a two-dimensional vector $\Phi=\Phi(\zeta,s)$, where the $2 \times 2$ matrices $A =
A(\zeta, s)$ and $B = B(\zeta,s)$ are given by
\begin{equation}\label{in:9}
  A(\zeta,s) = \left( 
    \begin{array}{c c}
      4 \zeta q &  4 \zeta^2 + s + 2q^2 + 2q'\\
      -4 \zeta^2 - s - 2 q^2 + 2q' & -4 \zeta q
    \end{array}\right)\qand
  B(\zeta,s) = \left( 
    \begin{array}{c c}
      q &  \zeta \\
      -\zeta & - q
    \end{array}\right).
\end{equation}
The compatibility of this overdetermined system implies that $q(s)$ solves Painlev\'e II.
Now let $\Phi =\begin{pmatrix} \Phi_1 \\ \Phi_2 \end{pmatrix}$ be the unique solution of
\eqref{in:8} satisfying
\begin{equation}\label{in:13}
  \Phi_1(\z;s)=\cos\left(\frac{4}{3}\z^3+s\z\right)+O(\z^{-1}), \quad \Phi_2(\z;s)=-\sin\left(\frac{4}{3}\z^3+s\z\right)+O(\z^{-1})\,,
\end{equation}
as $\z \to \pm \infty$ for $s\in \rr$.
The formula of \cite{schehr} is
\begin{equation}\label{eq:schehrgamma}
  \gamma(t,m)=\tfrac{16}{\pi^2} \langle
  h_{4^{2/3}t},h_{-4^{2/3}t}\rangle_{L^2(m,\infty)}
\end{equation}
where
\begin{equation}
  h_t(x)= \int_0^\infty d\zeta\,\zeta \Phi_2(\zeta,x) e^{-t\zeta^2}.
\end{equation}
Although Schehr's argument is non-rigourous, a later paper of \citet{baikLiechtySchehr}
proved directly the equivalence of the formula of \cite{schehr} and \eqref{eq:eff}, thus
establishing the validity of \eqref{eq:schehrgamma} based on Theorem \ref{thm:fPE}.

Before turning to the tail behavior of $\ct$, let us briefly explain how the uniqueness of
the maximizer of $\aip(t)-t^2$ can be established directly from the argument we described
above. In the derivation of the formula we assumed that the maximum of $\aip(t)-t^2$
is obtained at a unique point. However, it is not necessary to do this. In fact, if one
follows the argument without this assumption, one ends up with a formula for what is in
principle a super-probability density, i.e. a non-negative function $f(t,m)$ on
$\mathbb{R}\times\mathbb{R}$ with $\int_{\mathbb{R}\times\mathbb{R}}dm\,dt\, f(t,m)\ge 1$,
and in fact one can see from the argument that
\begin{equation}\label{eq:numberMax}
  \int_{\mathbb{R}\times\mathbb{R}}dm\,dt\, f(t,m) = {\rm expected~number~of ~maxima~of ~}\aip(t)-t^2.
\end{equation}
Recall that from \eqref{eq:johGOE} that the distribution of $\cm$ is given by a scaled
version of $F_{\rm GOE}$.  A non-trivial computation (see Section \pref{sec:uniq}) gives
\begin{equation}
  \int_{-\infty}^\infty dt\, f(t,m) = 4^{1/3}F'_\mathrm{GOE}(4^{1/3}m).
\end{equation}
This shows that $f(t,m)$ has total integral 1, which can only be true if the maximum is
unique almost surely, since the global maximum is attained at at least one point.

We mentioned earlier the conjecture that $\ct$ should have tails which decay like
$e^{-ct^3}$ (see e.g. \cite{halpZhang}). This can be proved using the techniques described
in this review:

\begin{thm}[\cite{quastelRemTails,schehr,baikLiechtySchehr,corwinHammond}]\label{thm:tail}
  There is a $c>0$ such that for every $\kappa>\frac{32}3$ and large enough $t$,
  \[e^{-\kappa t^3}\leq\pp\big(|\ct|>t\big)\leq
  ce^{-\frac43t^3+2t^2+\mathcal{O}(t^{3/2})}.\]
\end{thm}

\cite{corwinHammond} had obtained the $e^{-ct^3}$ decay for some $c>0$. The statement we
included here is the one appearing in \cite{quastelRemTails}. In fact, Schehr's formula
and its validation in \cite{baikLiechtySchehr} later yielded a lower bound that matches
the $e^{-\frac43t^3}$ behavior of the upper bound, so we know now that $\frac43$ is the
correct exponent. A precise asymptotic expansion of $\pp\big(|\ct|>t\big)$ based on that
formula has recently been obtained in \cite{bothnerLiechty}. The reason why
\cite{quastelRemTails} obtained a slightly worse lower bound is technical, and arises from
the fact that the explicit formula \eqref{eq:eff} for $f(t,m)$ is not useful for providing
a lower bound, and instead one needs to use a different argument.  On the other hand, the
upper bound can be obtained directly from \eqref{eq:eff}. In fact, the second formula
expresses this joint density as the difference of two Fredholm determinants, so we may use
\eqref{contintr} to estimate the difference, and then all that remains is to show that
this estimate can be integrated in $m$. See \cite{quastelRemTails} for more details.

\subsection{Local behavior of Airy\texorpdfstring{$_1$}{1}}
\label{sec:airy1}

As we mentioned, the boundary value kernel formulas introduced in Section \ref{sec:bdVlK} are better adapted
than the standard extended kernel formulas to study short range properties of the
processes. An interesting application is the following:

\begin{thm}[\cite{quastelRemAiry1}, Theorem \orefn{thm:regularity}]\label{thm:regularity}
  The Airy$\hspace{0.05em}_1$ process $\aipo$ and the Airy$\hspace{0.05em}_2$ process
  $\aip$ have versions with H\"older continuous paths with exponent $\tfrac12-\delta$ for
  any $\delta>0$.
\end{thm}

Continuity was known for $\aip$ (see Theorem \ref{thm:Airy2LPP}) but not for $\aipo$. The
H\"older $\frac12-$ continuity for $\aip$ also follows from the work of
\citet{corwinHammond}. Their proof is based on a certain Brownian Gibbs property for the
Airy$_2$ line ensemble (an infinite collection of continuous, non-intersecting paths, the
top line of which is $\aip$), and as such it cannot be extended to Airy$_1$, given that no
analog of the Airy$_2$ line ensemble is known in the flat case. This regularity is
expected to hold in fact for all the Airy processes in view of the fact that they are
believed to look locally like a Brownian motion (see Section \ref{sec:she}). Analogous
results have recently become available for the solutions of the KPZ equation at finite
times with certain initial conditions \cite{hairer,qrLocalBrownianKPZ,corwinHammondKPZ}.

The proof of Theorem \ref{thm:regularity} is based on an application of a suitable version
of the Kolmogorov criterion. In the Airy$_1$ case, it involves studying a truncated
version of the process, $\aipo^M(t)=\aipo(t)\uno{|\aipo(t)|\leq
  M}+M\uno{\aipo(t)>M}-M\uno{\aipo(t)<-M}$ and then proving the following estimate: for
fixed $\delta>0$, there is a $t_0\in(0,1)$ and an $n_0\in\nn$ such that for $0<t<t_0$,
  $n\geq n_0$ and $M=\big(3\log(t^{-(1+n)})\big)^{1/3}$ we have
 \[\ee\!\left(\big[\aipo^{M}(t)-\aipo^{M}(0)\big]^{2n}\right)\leq ct^{1+(1-\delta)n}\]
 where the constant $c>0$ is independent of $\delta$, $n_0$ and $t_0$. The proof of this
 estimate can be reduced to obtaining a suitable estimate on the difference
\[\left| \det(I-B_0+\bar P_ae^{t\Delta}\bar P_b e^{-t\Delta}B_0)-\det(I-B_0+\bar P_a
  B_0)\right|\] for $b\geq a\geq-M$. An important technical problem is that the kernels
appearing inside these determinants are not trace class, so one needs to conjugate
appropriately. We refer to \cite{quastelRemAiry1} for the details. The argument for
Airy$_2$ is similar.

As we mentioned, the Airy processes are expected to look locally like a Brownian
motion. In this direction, it can be shown using the boundary value kernel formulas that
the finite dimensional distributions of the Airy$_1$ process converge under
diffusive scaling to those of a Brownian motion. The same result was proved earlier by
\citet{hagg} for Airy$_2$ using different techinques. In fact, for Airy$_2$ a stronger
statement is now available (\citet{corwinHammond}), namely that it is locally absolutely
continuous with respect to Brownian motion.

\begin{thm}[\cite{quastelRemAiry1}, Theorem \orefn{thm:localBM}]\label{thm:localBM}
  For any fixed $s\in\mathbb{R}$, 
  let $B_\ep (\cdot)$ be defined by $B_\ep(t)= \ep^{-1/2}(\aipo(s+ \ep t)-\aipo(s) )$,
  $t>0$.  Then $B_\ep (\cdot)$ converges to Brownian motion in the sense of convergence of
  finite dimensional distributions.  The same holds for $\tilde{B}_\ep (\cdot)$ defined by
  $\tilde{B}_\ep(t)= B_\ep(-t)$, $t>0$.
\end{thm}

The proof of this result follows from an explicit computation of
\[\pp\!\left(\aipo(\ep t_1)\leq x+\sqrt{\ep}y_1 , \ldots, \aipo(\ep t_n)\leq
  x+\sqrt{\ep}y_n\,\middle|\,\aipo(0)=x\right)\]
and its limit as $\ep\to0$, see \cite{quastelRemAiry1} for the details. The same proof
works for the Airy$_2$ process and, in view of \eqref{eq:2to1}, it should be simple to
adapt it to the Airy$_{2\to1}$ process.

\subsection{Marginals of Airy\texorpdfstring{$_{2\to1}$}{2->1}}
\label{sec:airy2to1}

The last application of the results of Section \ref{sec:contStat} that we will discuss is
a proof of the conjecture \eqref{eq:conj} that the marginals of the Airy$_{2\to1}$ process
can be obtained from a variational problem for $\aip(t)-t^2$ on a half-line. The result is
the following:

\begin{thm}[\cite{qr-airy1to2}, Theorem 1]\label{thm:1to2}
  Fix $\alpha\in\rr$. For every $m\in\rr$,
  \[\pp\!\left(\sup_{t\leq\alpha}\big(\aip(t)-t^2\big)\leq m-\min\{0,\alpha\}^2\right)=\pp\!\left(\Bt(\alpha)\leq m\right).\]
\end{thm}

The right hand side can be expressed in terms of a Fredholm determinant. Define the
\emph{crossover distributions} $G^{2\to1}_\alpha$, for $\alpha\in\rr$, as
\[G^{2\to1}_\alpha(m)=\pp\!\left(\Bt(\alpha)\leq m\right).\]
We claim that
\begin{equation}
G^{2\to1}_\alpha(m)=\det\!\big(I-P_mK_\alpha P_m\big),\label{eq:Ga}
\end{equation}
where $K_\alpha=K^1_\alpha+K^2_\alpha$ and the kernels $K^1_\alpha$ and $K^2_\lambda$ are given
by
\[K_\alpha^1(x,y)=\int_0^{\infty} d\lambda\,e^{2\alpha\lambda}\Ai(x-\lambda+\max\{0,\alpha\}^2)\Ai(y+\lambda+\max\{0,\alpha\}^2)\]
and
\[K_\alpha^2(x,y)= \int_0^\infty
d\lambda\Ai(x+\lambda+\max\{0,\alpha\}^2)\Ai(y+\lambda+\max\{0,\alpha\}^2).\]
As noted in Appendix A of \cite{bfs}, the kernel $K_{2\to1}^{\rm ext}$ defined in \eqref{eqKCompleteInfinity} can be
expressed in terms of Airy functions:
\begin{equation}
  K_{2\to1}^{\rm ext}(s,t;x,y)=L_0(s,x;t,y)+e^{2t^3/3-2s^3/3+t\tilde y-s\tilde x}[L_1+L_2](s,x;t,y),\label{eq:Kinfty2}
\end{equation}
where
\begin{align}
      L_0(s,x;t,y)&=-e^{(s-t)\Delta}(\tilde x,\tilde
      y)=-\frac{1}{\sqrt{4\pi(t-s)}}e^{-(\tilde x-\tilde y)^2/4(t-s)},\\
      L_1(s,x;t,y)&=\int_0^\infty d\lambda\,e^{\lambda(s+t)}\Ai(\hat x-\lambda)\Ai(\hat
      y+\lambda),\\
      L_2(s,x,t,y)&=\int_0^\infty d\lambda\,e^{\lambda(t-s)}\Ai(\hat x+\lambda)\Ai(\hat
    y+\lambda)
\end{align}
with $\tilde x=x-s^2\uno{s\leq0}$, $\tilde y=y-t^2\uno{t\leq0}$, $\hat
x=x+s^2\uno{s\geq0}$ and $\hat y=y+t^2\uno{t\geq0}$.  Using this for $s=t=\alpha$ it is
straightforward to check that $K_{2\to1}^{\rm ext}(t,\cdot;t,\cdot)$ is just a conjugation
of the kernel $K_\alpha$, and \eqref{eq:Ga} follows.

The fact that $G^{2\to1}_\alpha$ crosses over between the GUE and GOE distributions is of
course a particular case of the crossover property of the Airy$_{2\to1}$ process, but can
be easily obtained from \eqref{eq:Ga} as well (see the discussion after Theorem 1 in
\cite{qr-airy1to2}).

The proof of Theorem \ref{thm:1to2} is similar to (and, in fact, somewhat simpler than)
the proof of \eqref{eq:Airy2GOE}. Basically, one applies Theorem \ref{thm:aiL} and the cyclic property
of determinants to compute the desired probability as
\begin{equation}
  \lim_{L\to\infty}\pp\!\left(\aip(t)\leq g(t)\text{ for }t\in[-L,\alpha]\right)
  =\lim_{L\to\infty}\det\!\left(I-K_{\Ai}+e^{(\alpha+L)H}\K\Theta^g_{[-L,\alpha]}\K\right)
  \label{eq:aiL}
\end{equation}
with $g(t)=t^2+\m$ and $\m=m-\min\{0,\alpha\}^2$. An argument similar to the one used in
Section \ref{sec:goe} (applying the Baker-Campbell-Hausdorff formula and later checking
the result rigorously, plus some asymptotic analysis to show that an error term goes to 0
in trace class norm as $L\to\infty$) yields
\begin{multline}\label{eq:limL}
  \pp\!\left(\sup_{t\leq\alpha}\big(\aip(t)-t^2\big)\leq m-\min\{0,\alpha\}^2\right)\\
  =\det\!\left(I-\K P_{\m+\alpha^2}K_{\Ai}-\K e^{\alpha\xi}\varrho_{\m+\alpha^2}
    e^{-\alpha\xi}\bar P_{\m+\alpha^2}\K\right).
\end{multline}
Since $\K=B_0P_0B_0$ and $B_0^2=I$ we have by the cyclic property of determinants that the right hand side of
\eqref{eq:limL} equals
\[\det\!\left(I-P_0B_0P_{\m+\alpha^2}B_0P_0-P_0B_0e^{\alpha\xi}
     \varrho_{\m+\alpha^2}e^{-\alpha\xi}\bar P_{\m+\alpha^2}B_0P_0\right).\]
Shifting the variables in the last determinant by $-m$ we deduce that
\begin{equation}
  \label{eq:modLim}
  \pp\!\left(\sup_{t\geq\alpha}\big(\aip(t)-t^2\big)\leq\m\right)
  =\det\!\left(I-P_mE_1P_m-P_mE_2P_m\right),
\end{equation}
where
\begin{align}
  E_1(x,y)&=\int_{-\infty}^{\m+\alpha^2}d\lambda
  \Ai(x-m+2\m+2\alpha^2-\lambda)e^{-2(\lambda-\m-\alpha^2)\alpha}\Ai(y-m+\lambda)\\
  \shortintertext{and}
  E_2(x,y)&=\int_{\m+\alpha^2}^\infty d\lambda \Ai(x-m+\lambda)\Ai(y-m+\lambda).
\end{align}
Shifting $\lambda$ by $\m+\alpha^2$ in both integrals and changing $\lambda$ to $-\lambda$
shows that $E_1(x,y)=K^1_\alpha(y,x)$ and $E_2=K^2_\alpha$,
whence the equality in Theorem \ref{thm:1to2} follows since $E_1^*=K^1_\alpha$ and $E_2^*=K_\alpha^2$.

\mbox{}

\printbibliography[heading=apa]

\end{document}